\documentclass[a4paper,twoside]{article}

\def\shorttitle{Inverse Source Problem for Hyperbolic Equations}
\def\shortauthor{D. Jiang, Y. Liu and M. Yamamoto}

\usepackage{amsfonts}
\usepackage{amsmath}
\usepackage{amssymb}
\usepackage{amscd}
\usepackage{amsbsy}
\usepackage{amsthm}
\usepackage{bm}
\usepackage{empheq}
\usepackage{graphicx}
\usepackage{psfrag}
\usepackage{color}
\usepackage{cite}
\usepackage{multicol}

\allowdisplaybreaks[4]

\newfont{\myfnt}{cmssi10 scaled 1440}
\numberwithin{equation}{section}
\catcode`@=11
\def\ps@nk{\def\@oddhead{\vbox{\hbox to \hsize{\pic \footnotesize \it \shorttitle
\hfill \rm \thepage} \vspace{1mm} \vspace*{-2mm}}}
\def\@evenhead{\vbox{\hbox to \hsize{\pic \footnotesize \rm \thepage \hfill \it \shortauthor}
\vspace{1mm} \vspace*{-2mm}}}
\def\@oddfoot{} \def\@evenfoot{}}
\def\ps@first{\def\@oddhead{\vbox{\hbox to \hsize{\pic \footnotesize
} \break}}
\def\@oddfoot{} \def\@evenfoot{}}
\newtheoremstyle{thmstyle}
  {6pt}
  {6pt}
  {\it}
  {}
  {\bf}
  {}
  {.5em}
  {}
\newtheoremstyle{remstyle}
  {6pt}
  {6pt}
  {\rm}
  {}
  {\bf}
  {}
  {.5em}
  {}
\def\Section#1{\Sec{\large #1} \setcounter{equation}{0} \vskip -6mm \indent}
\def\Sec{\@Startsection{section}{1}{\z@}
                                   {-3.5ex \@plus -1ex \@minus -.2ex}%
                                   {2.3ex \@plus.2ex}%
                                   {\normalfont\large\bfseries\boldmath}}
\def\@Startsection#1#2#3#4#5#6{%
  \if@noskipsec \leavevmode \fi
  \par
  \@tempskipa #4\relax
  \@afterindenttrue
  \ifdim \@tempskipa <\z@
    \@tempskipa -\@tempskipa \@afterindentfalse
  \fi
  \if@nobreak
    \everypar{}%
  \else
    \addpenalty\@secpenalty\addvspace\@tempskipa
  \fi
  \@ifstar
    {\@ssect{#3}{#4}{#5}{#6}}%
    {\@dblarg{\@Sect{#1}{#2}{#3}{#4}{#5}{#6}}}}
\def\@Sect#1#2#3#4#5#6[#7]#8{%
  \ifnum #2>\c@secnumdepth
    \let\@svsec\@empty
  \else
    \refstepcounter{#1}%
    \protected@edef\@svsec{\@seccntformat{#1}\relax}%
  \fi
  \@tempskipa #5\relax
  \ifdim \@tempskipa>\z@
    \begingroup
      #6{%
          \@hangfrom{\hskip #3\relax\@svsec \hskip -2.5mm}%
          \interlinepenalty \@M #8\@@par}
    \endgroup
    \csname #1mark\endcsname{#7}%
    \addcontentsline{toc}{#1}{%
      \ifnum #2>\c@secnumdepth \else
        \protect\numberline{\csname the#1\endcsname}%
      \fi
      #7}%
  \else
    \def\@svsechd{%
      #6{\hskip #3\relax
      \@svsec #8}%
      \csname #1mark\endcsname{#7}%
      \addcontentsline{toc}{#1}{%
        \ifnum #2>\c@secnumdepth \else
          \protect\numberline{\csname the#1\endcsname}%
        \fi
        #7}}%
  \fi
  \@xsect{#5}}
\renewenvironment{abstract}{%
        \small
        \quotation
         \noindent {\bfseries \abstractname } }%
      {\if@twocolumn\else\endquotation\fi}
\def\Subsection#1{\Subsec{#1} \vskip -6mm \indent}
\def\Subsec{\@StartSubsection{subsection}{2}{\z@}%
                                     {-3.25ex\@plus -1ex \@minus -.2ex}%
                                     {1.5ex \@plus .2ex}%
                                     {\normalfont\normalsize\bfseries\boldmath}}
\def\@StartSubsection#1#2#3#4#5#6{%
  \if@noskipsec \leavevmode \fi
  \par
  \@tempskipa #4\relax
  \@afterindenttrue
  \ifdim \@tempskipa <\z@
    \@tempskipa -\@tempskipa \@afterindentfalse
  \fi
  \if@nobreak
    \everypar{}%
  \else
    \addpenalty\@secpenalty\addvspace\@tempskipa
  \fi
  \@ifstar
    {\@ssect{#3}{#4}{#5}{#6}}%
    {\@dblarg{\@SubSect{#1}{#2}{#3}{#4}{#5}{#6}}}}
\def\@SubSect#1#2#3#4#5#6[#7]#8{%
  \ifnum #2>\c@secnumdepth
    \let\@svsec\@empty
  \else
    \refstepcounter{#1}%
    \protected@edef\@svsec{\@seccntformat{#1}\relax}%
  \fi
  \@tempskipa #5\relax
  \ifdim \@tempskipa>\z@
    \begingroup
      #6{%
          \@hangfrom{\hskip #3\relax\@svsec\hskip -1.5mm}%
          \interlinepenalty \@M #8\@@par}
    \endgroup
    \csname #1mark\endcsname{#7}%
    \addcontentsline{toc}{#1}{%
      \ifnum #2>\c@secnumdepth \else
        \protect\numberline{\csname the#1\endcsname}%
      \fi
      #7}%
  \else
    \def\@svsechd{%
      #6{\hskip #3\relax
      \@svsec #8}%
      \csname #1mark\endcsname{#7}%
      \addcontentsline{toc}{#1}{%
        \ifnum #2>\c@secnumdepth \else
          \protect\numberline{\csname the#1\endcsname}%
        \fi
        #7}}%
  \fi
  \@xsect{#5}}
\def\list#1#2{\ifnum \@listdepth >5\relax \@toodeep \else \global
\advance \@listdepth\@ne \fi \rightmargin \z@ \listparindent\z@
\itemindent\z@ \csname @list\romannumeral\the\@listdepth\endcsname
\def\@itemlabel{#1}\let\makelabel\@mklab \@nmbrlistfalse #2\relax
\@trivlist \parskip 0pt \parindent\listparindent \advance \linewidth
-\rightmargin \advance\linewidth -\leftmargin \advance\@totalleftmargin
\leftmargin \parshape \@ne \@totalleftmargin \linewidth \ignorespaces}
\renewcommand{\@makecaption}[2]{\begin{center}#1. #2\end{center}}
\catcode`@=12 
\theoremstyle{thmstyle}
\newtheorem{thm}{\indent Theorem}[section]
\newtheorem{lem}[thm]{\indent Lemma}
\newtheorem{prop}[thm]{\indent Proposition}

\newtheorem{defi}[thm]{\indent Definition}

\newtheorem{prob}[thm]{\indent Problem}
\theoremstyle{remstyle}

\newtheorem{algo}[thm]{\indent Algorithm}
\newtheorem{ex}[thm]{\indent Example}

\newsavebox{\mygraphic}
\def\pic{\begin{picture}(0,0) \put(-210,-1250){\usebox{\mygraphic}} \end{picture}}
\newfont{\HUGEbf}{cmbx10 scaled 3500}
\definecolor{gray}{rgb}{0.9,0.9,0.9}
\def\thebibliography#1{\section*{\bf \large References}
\list{[\arabic{enumi}]} {\settowidth \labelwidth{[#1]} \leftmargin
\labelwidth \advance \leftmargin \labelsep \usecounter{enumi}}
\def\newblock{\hskip .11em plus .33em minus .07em} \footnotesize \sloppy \clubpenalty
4000 \widowpenalty 4000 \sfcode`\.=1000 \relax}

\def\BR{\mathbb R}
\def\cA{\mathcal A}
\def\cD{\mathcal D}
\def\cH{\mathcal H}
\def\a{\mathrm a}
\def\rb{\mathrm b}
\def\rc{\mathrm c}
\def\rd{\mathrm d}
\def\diam{\mathrm{diam}}
\def\rdiv{\mathrm{div}}
\def\e{\mathrm e}
\def\supp{\mathrm{supp}}
\def\true{\mathrm{true}}
\def\Ga{\Gamma}

\def\Om{\Omega}
\def\al{\alpha}
\def\be{\beta}
\def\ga{\gamma}
\def\de{\delta}
\def\ep{\epsilon}
\def\ve{\varepsilon}
\def\te{\theta}
\def\ka{\kappa}
\def\la{\lambda}
\def\vp{\varphi}
\def\om{\omega}
\def\f{\frac}
\def\nb{\nabla}
\def\ov{\overline}
\def\pa{\partial}
\def\tri{\triangle}
\def\wt{\widetilde}

\theoremstyle{definition}

\numberwithin{equation}{section}

\title{\Large\bf\boldmath Theoretical Stability and Numerical Reconstruction for an Inverse Source Problem for Hyperbolic Equations$^*$}

\author{\large Daijun JIANG$^\dag$\qquad Yikan LIU$^\ddag$\qquad Masahiro YAMAMOTO$^\ddag$}

\date{}

\begin{document}

\maketitle

\thispagestyle{first}
\renewcommand{\thefootnote}{\fnsymbol{footnote}}

\footnotetext{\hspace*{-5mm} \begin{tabular}{@{}r@{}p{14cm}@{}} &
Manuscript last updated: \today.\\
$^\dag$ & Mathematics and Statistics, Central China Normal University, Wuhan 430079, China.\\
& E-mail: jiangdaijun@mail.ccnu.edu.cn.\\
$^\ddag$ & Graduate School of Mathematical Sciences, The University of Tokyo, 3-8-1 Komaba, Meguro-ku, Tokyo 153-8914, Japan. E-mail: ykliu@ms.u-tokyo.ac.jp, myama@ms.u-tokyo.ac.jp\\
$^*$ & This work has been partly supported by the A3 Foresight Program ``Modeling and Computation of Applied Inverse Problems'', Japan Society of the Promotion of Science (JSPS). D.\! Jiang has been financially supported by self-determined research funds of CCNU from the colleges' basic research and operation of MOE (No.\! CCNU14A05039), National Natural Science Foundation of China (Nos.\! 11326233, 11401241 and 11571265) and China Postdoctoral Science Foundation (Grant No.\! 2012M521444). Y.\! Liu and M.\! Yamamoto have been supported by Grant-in-Aid for Scientific Research (S) 15H05740, JSPS.
\end{tabular}}

\renewcommand{\thefootnote}{\arabic{footnote}}

\begin{abstract}
In this paper, we investigate the inverse problem on determining the spatial component of the source term in a hyperbolic equation with time-dependent principal part. Based on a newly established Carleman estimate for general hyperbolic operators, we prove a local stability result of H\"older type in both cases of partial boundary and interior observation data. Numerically, we adopt the classical Tikhonov regularization to transform the inverse problem into an output least-squares minimization, which can be solved by the iterative thresholding algorithm. The proposed algorithm is computationally easy and efficient: the minimizer at each step has explicit solution. Abundant amounts of numerical experiments are presented to demonstrate the accuracy and efficiency of the algorithm.

\vskip 4.5mm

\noindent\begin{tabular}{@{}l@{ }p{10cm}} {\bf Keywords } & Inverse source problem, Hyperbolic equation, Carleman estimate,\\
& Iterative thresholding algorithm
\end{tabular}

\end{abstract}

\baselineskip 14pt

\setlength{\parindent}{1.5em}

\setcounter{section}{0}

\Section{Introduction}\label{sec-intro}

Let $\Om\subset\BR^n$ ($n=1,2,\ldots$) be an open bounded domain with a smooth boundary $\pa\Om$ (e.g., of $C^2$-class), and let $\nu=\nu(x)=(\nu_1(x),\ldots,\nu_n(x))$ be the outward unit normal vector to $\pa\Om$ at $x\in\pa\Om$. For some $T>0$, set $Q:=\Om\times(-T,T)$. We consider the following initial value problem for a hyperbolic equation whose principal part depends on the time variable
\begin{empheq}[left=\empheqlbrace,right=]{alignat=2}
& (\pa_t^2-\cA(t))u(x,t)=F(x,t) & \quad & ((x,t)\in Q),\label{eq-gov-u}\\
& u(x,0)=\pa_tu(x,0)=0 & \quad & (x\in\Om),\label{eq-IC-u}
\end{empheq}
where
\begin{align*}
\cA(t)u(x,t) & :=\rdiv(a(x,t)\nb u(x,t))+b(x,t)\cdot\nb u(x,t)+c(x,t)u(x,t)\\
& \:=\sum_{i,j=1}^n\pa_j(a_{ij}(x,t)\pa_ku(x,t))+\sum_{i=1}^nb_i(x,t)\pa_ju(x,t)+c(x,t)u(x,t).
\end{align*}
Here $a=(a_{ij})_{1\le i,j\le n}$ is a symmetric matrix, $b=(b_i)_{1\le i\le n}$ is a vector, and there exists a constant $\ka>0$ such that
\[a(x,t)\xi\cdot\xi=\sum_{i,j=1}^na_{ij}(x,t)\xi_i\xi_j\ge\ka|\xi|^2=\ka\sum_{i=1}^n\xi_i^2\quad(\forall\,(x,t)\in\ov Q\,,\ \forall\,\xi=(\xi_1,\ldots,\xi_n)\in\BR^n).\]
The regularities of $a,b,c$, the assumptions on the source term $F$ and the boundary condition will be specified later. We denote the normal derivative associated with the elliptic operator $\cA(t)$ as
\[\pa_\cA u:=a\nu\cdot\nb u\quad\mbox{on }\pa\Om\times(-T,T).\]
The well-posedness result concerning \eqref{eq-gov-u} will be provided in Lemma \ref{lem-Lions}.

The main focuses of this paper are the theoretical stability and the numerical treatment for the following inverse source problem.

\begin{prob}\label{prob-ISP}
Let the subboundary $\Ga\subset\pa\Om$, the subdomain $\om\subset\Om$ and $T>0$ be suitably given. Assume that the source term $F(x,t)=f(x)R(x,t)$ in \eqref{eq-gov-u} where $R$ is given, and let $u$ satisfy \eqref{eq-gov-u}--\eqref{eq-IC-u}. Determine $f(x)$ by

{\bf Case (I).} the partial boundary observation data $\{u,\pa_\cA u\}|_{\Ga\times(-T,T)}$, or

{\bf Case (II).} the partial interior observation data $u|_{\om\times(-T,T)}$.
\end{prob}

Investigating the above problem from both theoretical and numerical aspects not only originates from the interest on mathematics, but also roots in its significance in practice. In our formulation \eqref{eq-gov-u}, the source term $f(x)R(x,t)$ is incompletely separated into its spatial and temporal components, and the purposed inverse problem means the determination of the spatial component $f$. Especially, if the source term is in form of complete separation of variables, i.e.\! $R$ is space-independent, \eqref{eq-gov-u} becomes an approximation to a model for elastic waves, and the term $f(x)R(t)$ acts as the external force modeling vibrations (see Yamamoto \cite{Y95}). Recently, it reveals in Liu and Yamamoto \cite{LY14} that the one-dimensional time cone model for crystallization growth (see Cahn \cite{C95}) indeed takes the form of \eqref{eq-gov-u}, where the principal part involves the time-dependent growth speed, and $f(x)$ stands for the spatial distribution of the nucleation rate.

Although inverse hyperbolic problems have attracted considerable attentions during the last two decades, the majority of existing works only treated inverse source problems like Problem \ref{prob-ISP} with time-independent principal part, which is technically easier to show the uniqueness and the stability by Carleman estimates. We refer to Bukhgeim and Klibanov \cite{BK81} for the uniqueness, and Yamamoto \cite{Y95,Y99}, Puel and Yamamoto \cite{PY96}, Isakov and Yamamoto \cite{IY00}, Imanuvilov and Yamamoto \cite{IY01a,IY01b} for the stability. Especially, the global Lipschitz stability was proved for the boundary measurement case in \cite{IY01b} and the interior measurement case in \cite{IY01a} under the probably optimal geometrical condition on observable regions. However, if the principal part is time-dependent, we need to argue extra and there are no publications to the best knowledge of the authors. We can refer to Li and Yamamoto \cite{LY13} as a related work. The present paper is mainly motivated by \cite{IY01a,IY01b} to establish similar stability results under the more general setting on the principal part by a refined Carleman estimate. For comprehensive discussions on inverse hyperbolic problems by Carleman estimates, see Bellassoued and Yamamoto \cite{BY14}.

Correspondingly, works on numerical reconstructions of source terms in hyperbolic equations are quite limited compared with those of coefficients. Regarding the numerical approaches to coefficient inverse hyperbolic problems and related topics, we refer to the two monographs \cite{KT04,KSS04}. In \cite{LXY12}, the authors developed a spectral method for the inverse coefficient problem for the hyperbolic equation derived in \cite{LY14}. On the other hand, a class of iterative thresholding algorithms was purposed for linear inverse problems in early 2000s, whose convergence was first rigorously analyzed in Daubechies, Defrise and De Mol \cite{DDM04}. As an extension of classical gradient algorithms with regularization, the iterative thresholding algorithm and its updated versions have proved their feasibility mainly in the abundant applications to image processing due to their simplicity (see \cite{BT09,BF07,DTV07,FR08}). However, the flavor of this method is less familiar among the researchers of inverse problems for partial differential equations. Recently, the iterative thresholding algorithm was utilized in Jiang, Feng and Zou \cite{JFZ14} to treat inverse problems for elliptic and parabolic equations. Very recently, based on the theoretical stability of Lipschitz type, we develop a similar iterative method for an inverse source problem in the three-dimensional time cone model in \cite{LJY15}.

In this paper, we first establish a new Carleman estimate for the hyperbolic operator $\pa_t^2-\cA(t)$ in \eqref{eq-gov-u} which estimates also second-order derivatives, by which we can prove the local H\"older stability for both boundary and interior observations in accordance with the observation time $T$ and the observable region $\Ga$ or $\om$. Numerically, by interpreting Problem \ref{prob-ISP} as a minimization problem, we characterize the minimizer by a variational equation which involves the adjoint problem of \eqref{eq-gov-u}. This leads to the desired iterative thresholding algorithm, and we test its performance from various aspects by numerical experiments up to three spatial dimensions.

The remainder of this paper is organized as follows. In Section \ref{sec-stab}, we recall some existing results related to \eqref{eq-gov-u} and state the main stability result of H\"older type concerning Problem \ref{prob-ISP}. Then Section \ref{sec-proof} is devoted to the proof of the main result. In Section \ref{sec-recon}, we reformulate our inverse source problem as an optimization problem for numerical treatments, and purpose the iteration thresholding algorithm. Abundant numerical tests along with discussions on the performance are carried out in Section \ref{sec-numer}, and concluding remarks will be given in Section \ref{sec-con}. Finally, the key Carleman estimates used in Section \ref{sec-proof} is proved in Appendix \ref{sec-app}.

\Section{Preliminary and Theoretical Stability}\label{sec-stab}

In this section, we first review some existing results on the forward and inverse problems related to \eqref{eq-gov-u}, and then give the statement of the main result on the theoretical stability.

Recall $Q:=\Om\times(-T,T)$. Let $H^k(Q)$, $H^{k-\f12}(\pa\Om)$, $W^{k,\infty}(\Om)$, etc.\! ($k=0,1,\ldots$) denote usual Sobolev spaces. Concerning the solvability and stability issues for \eqref{eq-gov-u}, as a typical example we state a well-posedness result of the following initial-boundary value problem
\begin{equation}\label{eq-ibvp-u}
\begin{cases}
(\pa_t^2-\cA(t))u=\pa_t^2u-\rdiv(a\nb u)-b\cdot\nb u-cu=F & \mbox{in }Q,\\
u=g_0,\ \pa_tu=g_1 & \mbox{in }\Om\times\{0\},\\
\pa_\cA u=h & \mbox{on }\pa\Om\times(-T,T).
\end{cases}
\end{equation}

\begin{lem}[see Isakov \cite{I06}]\label{lem-Lions}
{\rm(a)} Let $k=1,2,\ldots$. Assume that $\pa\Om\in C^{k+1}$, $a\in W^{k,\infty}(Q)$, $b,c\in W^{k-1,\infty}(Q)$, and
\[F\in H^{k-1}(Q),\quad g_0\in H^k(\Om),\quad g_1\in H^{k-1}(\Om),\quad h\in H^{k-\f12}(\pa\Om\times(-T,T))\]
satisfy the $k$th-order compatibility condition on $\pa\Om\times\{0\}$. Then there is a unique solution $u$ to the initial-boundary value problem \eqref{eq-ibvp-u}. Moreover, there exists a constant $C=C(\cA(t),\Om,T)>0$ such that
\begin{align}
\|u(\,\cdot\,,t)\|_{H^k(\Om)}+\|\pa_tu(\,\cdot\,,t)\|_{H^{k-1}(\Om)} & \le C\Big(\|F\|_{H^{k-1}(Q)}+\|g_0\|_{H^k(\Om)}+\|g_1\|_{H^{k-1}(\Om)}\nonumber\\
& \qquad\:\:\:+\|h\|_{H^{k-\f12}(\pa\Om\times(-T,T))}\Big).\label{eq-est-u}
\end{align}

{\rm(b)} Moreover, in the case of $k=1$, assume that $\pa_ta=0$, $\pa_tb,\pa_tc\in L^\infty(Q)$, and
\[\pa_tF\in L^2(Q),\quad g_0\in H^2(\Om),\quad g_1\in L^2(\Om),\quad\pa_th\in H^{\f12}(\pa\Om\times(-T,T))\]
satisfy the second-order compatibility condition, then $\|\pa_t^2u(\,\cdot\,,t)\|_{L^2(\Om)}$ is bounded by the right-hand side of \eqref{eq-est-u} and by
\[\|\pa_tF\|_{L^2(Q)}+\|g_0\|_{H^2(\Om)}+\|g_1\|_{L^2(\Om)}+\|\pa_th\|_{H^{\f12}(\pa\Om\times(-T,T))}.\]
\end{lem}

The above lemma follows in principle from Lions and Magenes \cite{LM72}, and similar results also holds with other types of inhomogeneous boundary conditions. However, the regularity assumptions in Lemma \ref{lem-Lions} are merely a sufficient condition guaranteeing the well-posedness because the optimal one is still unknown.

Regarding the related inverse source problem (i.e., Problem \ref{prob-ISP}), the global Lipschitz stability is well-known for the wave equation. As a representative, we consider
\begin{equation}\label{eq-wave}
\begin{cases}
\pa_t^2u(x,t)=\tri u(x,t)+p(x)u(x,t)+f(x)R(x,t) & (x\in\Om,\ 0<t<T),\\
u(x,0)=\pa_tu(x,0)=0 & (x\in\Om),\\
\pa_\nu u(x,t):=\nb u(x,t)\cdot\nu(x)=0 & (x\in\pa\Om,\ 0<t<T)
\end{cases}
\end{equation}
and state the conclusion as follows.

\begin{lem}[see Imanuvilov and Yamamoto \cite{IY01a}]\label{lem-inter}
Let $\om$ be a subdomain of $\Om$ such that
\begin{equation}\label{eq-cond}
\{x\in\pa\Om;\,(x-x_0)\cdot\nu(x)\ge0\}\subset\pa\om\mbox{ for some }x_0\notin\ov{\Om\setminus\om}\,,\mbox{ and }T>\sup_{x\in\Om}|x-x_0|.
\end{equation}
Further assume that for any $x\in\pa\Om\setminus\pa\om$, there exists an open ball $U_x$ centered at $x$ such that $U_x\cap\Om$ is convex. Let $u$ be the solution to the initial-boundary value problem \eqref{eq-wave}, where $f\in L^2(\Om)$, $p\in L^\infty(\Om)$, $R\in H^1(0,T;L^\infty(\Om))$ and there exist constants $M>0$, $r_0>0$ such that
\[\|p\|_{L^\infty(\Om)}\le M,\quad\|\pa_tR\|_{L^2(0,T;L^\infty(\Om))}\le M,\quad|R(\,\cdot\,,0)|\ge r_0\mbox{ in }\Om.\]
Then there exists a constant $C=C(M,r_0,\om,\Om,T)>0$ such that
\[\|f\|_{L^2(\Om)}\le C\left(\|\pa_tu\|_{L^2(\om\times(0,T))}+\|\pa_t^2u\|_{L^2(\om\times(0,T))}\right).\]
\end{lem}

\begin{figure}[htbp]\centering
\includegraphics[trim=10mm 8mm 8mm 10mm,clip=true,width=.45\textwidth]{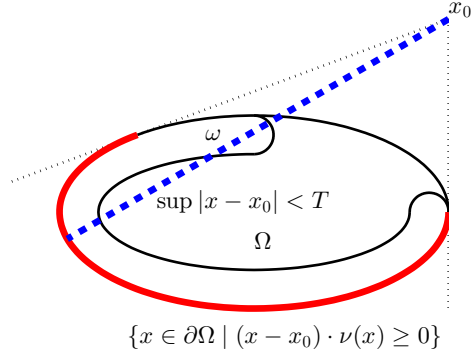}\\
\caption{A typical example for the spatial and temporal assumption guaranteeing the global Lipschitz stability of Problem \ref{prob-ISP} in case of the wave equation.}\label{fig-cond}
\end{figure}

The formulation \eqref{eq-wave} corresponds with the more general one of \eqref{eq-gov-u}--\eqref{eq-IC-u} in view of an even extension $u(x,-t):=u(x,t)$ for $t>0$. Lemma \ref{lem-inter} motivates the investigation on possible stability results for the generalized problem \eqref{eq-gov-u}--\eqref{eq-IC-u}. Meanwhile, it is also the starting point for developing numerical methods for the reconstruction of the source term, because the theoretical stability is guaranteed under condition \eqref{eq-cond} on the observable subdomain $\om$ and duration $T$. Although such a condition mainly originates from its necessity in the proof of Lemma \ref{lem-inter} by Carleman estimates, it also follows naturally from the essence of wave propagation. On the one hand, $\om$ cannot be too localized to capture waves in all directions. On the other hand, due to the finite propagation speed, adequate observation time should be given for the distant wave to reach $\om$. We illustrate a typical choice of $x_0$, $\om$ and $T$ in Figure \ref{fig-cond} for readers' better understanding.

Now we turn to the more general case of time-dependent principal part $\cA(t)$ in the governing equation \eqref{eq-gov-u}. We assume that the hyperbolic operator $\pa_t^2-\cA(t)$ admits a Carleman estimate. More precisely, for given $d\in C^2(\ov\Om)$ satisfying $d>0$ on $\ov\Om$, we set
\begin{equation}\label{eq-weight}
\psi(x,t):=d(x)-\be\,t^2,\quad\vp(x,t):=\e^{\la\psi(x,t)}
\end{equation}
with $\be\in(0,1)$ and a large parameter $\la>0$. Furthermore for $\de>0$, we set
\begin{equation}\label{eq-def-Qde}
Q_\de:=\{(x,t)\in Q;\,\psi(x,t)>\de\},\quad\Om_\de:=\{x\in\Om;\,\psi(x,0)>\de\}.
\end{equation}
Let a subboundary $\Ga\subset\pa\Om$ and a subdomain $\om\subset\Om$ satisfy
\begin{equation}\label{eq-asp-Geo}
\ov{Q_0}\subset(\Om\cup\Ga)\times(-T,T),\quad\pa\om\supset\Ga.
\end{equation}
We assume that a Carleman estimate holds with the weight function $\vp(x,t)$, that is, there exist constants $s_0>0$ and $C>0$ such that
\begin{equation}\label{eq-CE-0}
\int_{Q_0}s\left(|\pa_tu|^2+|\nb u|^2+s^2|u|^2\right)\e^{2s\vp}\,\rd x\rd t\le C\int_{Q_0}\left|(\pa_t^2-\cA(t))u\right|^2\e^{2s\vp}\,\rd x\rd t
\end{equation}
for all $s\ge s_0$ and all $u\in H^2(Q)$ satisfying $\supp\,u\subset Q_0$.

The Carleman estimate \eqref{eq-CE-0} definitely relies on some conditions on $d$ and the principal coefficients $a_{jk}$ (see Amirov and Yamamoto \cite{AY05}, Isakov \cite{I06}, Kha\u\i darov \cite{K87}).

Regarding Problem \ref{prob-ISP}, our main result is summarized as the following local H\"older estimate.

\begin{thm}\label{thm-ISP}
Assume that the hyperbolic operator $\pa_t^2-\cA(t)$ in \eqref{eq-gov-u} admits the Carleman estimate \eqref{eq-CE-0} with the weight function defined in \eqref{eq-weight}, and the coefficients of $\cA(t)$ satisfy
\[a\in W^{2,\infty}(-T,T;W^{1,\infty}(\Om)),\quad b,c\in W^{2,\infty}(-T,T;L^\infty(\Om)).\]
Let $u$ satisfy \eqref{eq-gov-u}--\eqref{eq-IC-u}, where $F(x,t)=f(x)R(x,t)$, $f\in L^2(\Om)$, $R\in W^{2,\infty}(-T,T;L^\infty(\Om))$ and there exists a constant $r_0>0$ such that
\begin{equation}\label{eq-asp-R}
|R(\,\cdot\,,0)|\ge r_0\quad\mbox{in }\Om\,.
\end{equation}
Further assume that
\begin{equation}\label{eq-reg-u}
u\in\bigcap_{k=0}^2H^{4-k}(-T,T;H^k(\Om))
\end{equation}
and there exists a constant $M>0$ such that
\begin{equation}\label{eq-pri-bdd}
\sum_{k=0}^2\|u\|_{H^{3-k}(-T,T;H^k(\Om))}\le M.
\end{equation}
Provided that $\Ga\subset\pa\Om$ and $\om\subset\Om$ satisfy condition \eqref{eq-asp-Geo}, then there exists $T_0>0$ satisfying: for arbitrarily given $\de>0$ and $T>T_0$, there exist constants $C>0$ and $\te\in(0,1)$ such that
\begin{equation}\label{eq-thm-est}
\|f\|_{L^2(\Om_\de)}\le CD+CM^{1-\te}D^\te,
\end{equation}
where $\Om_\de$ was defined in \eqref{eq-def-Qde} and
\begin{equation}\label{eq-def-obs}
D:=\begin{cases}
\|u\|_{H^4(-T,T;H^{-\f12}(\Ga))}+\|u\|_{H^2(-T,T;H^{\f32}(\Ga))}+\|\pa_\cA u\|_{H^2(-T,T;H^{\f12}(\Ga))}, & \mbox{Case }{\rm(I)},\\
\|u\|_{H^4(-T,T;L^2(\om))}+\|u\|_{H^2(-T,T;H^2(\om))}, & \mbox{Case }{\rm(II)}.
\end{cases}
\end{equation}
Here $T_0$ depends on $\Ga,\cA(t)$ and $C,\te$ depend on $T,\Ga,\cA(t)$ in Case {\rm(I)}, while $T_0$ depends on $\om,\cA(t)$ and $C,\te$ depend on $T,\om,\cA(t)$ in Case {\rm(II)}.
\end{thm}

In the above theorem, we assume comparatively low regularities of the coefficients in $\cA(t)$ and the source term $F$, which may not guarantee \eqref{eq-reg-u} in view of the forward problem. Actually, a sufficient condition for \eqref{eq-reg-u} is
\[\pa\Om\in C^5,\quad a\in W^{4,\infty}(Q),\quad b,c\in W^{3,\infty}(Q),\quad F\in H^3(Q)\]
by Lemma \ref{lem-Lions}(a) with $k=4$ and investigating the governing equation of $\pa_tu$ and $\pa_t^2u$. Further, if the principal part of $\cA(t)$ is independent of time, then the condition can be weakened to
\[\pa\Om\in C^2,\quad a\in W^{1,\infty}(\Om),\quad b,c\in W^{3,\infty}(-T,T;L^\infty(\Om)),\quad F\in\bigcap_{k=0}^2H^{3-k}(-T,T;H^k(\Om))\]
by Lemma \ref{lem-Lions}(b) and investigating the governing equations of $\pa_t^\ell u$ ($\ell=1,2,3$). However, in Theorem \ref{thm-ISP} it is understood that we are in advance given the observation data with certain smoothness in $\om\times(-T,T)$ or on $\Ga\times(-T,T)$, which does not necessarily solve the initial-boundary value problem \eqref{eq-ibvp-u} in $Q$. Meanwhile, we only assume the minimum necessary conditions on $a,b,c$ and $F$ which will be used in the proof of Theorem \ref{thm-ISP}.

For the parabolic case with $t$-dependent principal part, a similar stability was already proved in Imanuvilov and Yamamoto \cite{IY98}.  In fact, in the parabolic case, the proof for the $t$-dependent case is the same as the $t$-independent case because of the character of the parabolic Carleman estimate.

The condition \eqref{eq-asp-R} is essential. We cannot expect to determine the $(x,t)$-dependent $F(x,t)$ of \eqref{eq-gov-u} in general by our data of the inverse problem, because our extra measurement data $\pa_\cA u|_{\Ga\times(-T,T)}$ for the inverse problem depend on $n$ independent variables on $\Ga\times(-T,T)$ (i.e., $(n-1)$-spatial variables and one time derivative), but our unknown $F$ depends on $(n+1)$-independent variables $(x,t)$. Here we note that $u|_{\Ga\times(-T,T)}$, $u(\,\cdot\,,0)$ and $\pa_tu(\,\cdot\,,0)$ should be regarded as one boundary condition and the initial condition as the direct problem and so for determining a source term, we can consider only $\pa_\cA u|_{\Ga\times(-T,T)}$ as extra data.

Henceforth, by $C>0$ we denote a generic constant independent of the large parameter $s>0$ which may change from line to line, unless specified otherwise.

\Section{Proof of Theorem \ref{thm-ISP}}\label{sec-proof}

This section is devoted to the proof of Theorem \ref{thm-ISP} under the same assumptions therein. The argument is mainly based on the following key Carleman estimates for $\pa_t^2-\cA(t)$ in \eqref{eq-gov-u} which estimates also second-order derivatives.

\begin{lem}\label{lem-Carleman}
Assume that the hyperbolic operator $\pa_t^2-\cA(t)$ admits the Carleman estimate \eqref{eq-CE-0} with the weight function $\vp(x,t)$ in \eqref{eq-weight}, where the coefficients in $\cA(t)$ satisfy
\[a\in W^{1,\infty}(-T,T;W^{1,\infty}(\Om)),\quad b,c\in W^{1,\infty}(-T,T;L^\infty(\Om)).\]
Recall the set $Q_\de\subset Q=\Om\times(-T,T)$ defined in \eqref{eq-def-Qde}. Then for arbitrarily given $\de\ge0$, there exist constants $C>0$ and $s_0>0$ such that
\begin{align}
& \int_{Q_\de}\sum_{i,j=1}^n|\pa_i\pa_ju|^2\,\e^{2s\vp}\,\rd x\rd t\le C\int_{Q_\de}\left(\f1s|\pa_tF|^2+s\,|F|^2\right)\e^{2s\vp}\,\rd x\rd t\label{eq-CE-1}\\
& \int_{Q_\de}|\pa_t^2u|^2\,\e^{2s\vp}\,\rd x\rd t\le C\int_{Q_\de}\left(\f1s|\pa_tF|^2+|F|^2\right)\e^{2s\vp}\,\rd x\rd t,\label{eq-CE-2}
\end{align}
hold for all $s\ge s_0$ and all $u\in H^2(Q_\de)$ satisfying $\pa_tu\in H^2(Q_\de)$, $\pa_t^\ell u(\,\cdot\,,\pm T)=0$ for $\ell=0,1,2$ and $\supp\,u\subset Q_\de$, where $F:=\pa_t^2u-\cA(t)u$ in $Q$.
\end{lem}

For the sake of consistency, we postpone the proof of Lemma \ref{lem-Carleman} to Appendix \ref{sec-app}.

\begin{proof}[Completion of the proof of Theorem $\ref{thm-ISP}$]
By the Sobolev extension theorem, there exists $\wt u\in H^4(-T,T;L^2(\Om))\cap H^2(-T,T;H^2(\Om))$ such that
\[\begin{cases}
\wt u=u,\ \pa_\cA\wt u=\pa_\cA u\quad\mbox{on }\Ga\times(-T,T), & \mbox{Case (I)},\\
\wt u=u\quad\mbox{in }\om\times(-T,T), & \mbox{Case (II)}
\end{cases}\]
and
\begin{equation}\label{eq-est-wtu}
\sum_{k=0}^2\|\wt u\|_{H^{3-k}(-T,T;H^k(\Om))}\le CM,\quad\|\wt u\|_{H^4(-T,T;L^2(\Om))}+\|\wt u\|_{H^2(-T,T;H^2(\Om))}\le CD,
\end{equation}
where $M$ and $D$ were defined in \eqref{eq-pri-bdd} and \eqref{eq-def-obs}, respectively.

Recall $\cA(t)w=\rdiv(a\nb w)+b\cdot\nb w+cw$ and introduce
\begin{equation}\label{eq-def-A'}
\begin{aligned}
\cA'(t)w & :=\rdiv((\pa_ta)\nb w)+(\pa_tb)\cdot\nb w+(\pa_tc)w,\\
\cA''(t)w & :=\rdiv((\pa_t^2a)\nb w)+(\pa_t^2b)\cdot\nb w+(\pa_t^2c)w
\end{aligned}
\end{equation}
for later use. Setting
\[\cH:=\pa_t^2-\cA(t),\quad y:=u-\wt u,\quad y_1:=\pa_ty,\quad y_2:=\pa_t^2y,\]
we see $y\in\bigcap_{k=0}^2H^{4-k}(-T,T;H^k(\Om))$, and direct calculations yields
\begin{equation}\label{eq-gov-y}
\begin{aligned}
\cH y & =fR-\cH\wt u,\quad\cH y_1=\cA'(t)y+f\,\pa_tR-\pa_t(\cH\wt u\,),\\
\cH y_2 & =2\cA'(t)y_1+\cA''(t)y+f\,\pa_t^2R-\pa_t^2(\cH\wt u\,).
\end{aligned}
\end{equation}

Since it suffices to consider sufficiently small $\de>0$, by $d>0$ in $\ov\Om\,$, we can choose $\de>0$ such that
\begin{equation}\label{eq-asp-de}
\min_{x\in \ov\Om} d(x) \ge 2\de.
\end{equation}
We define a cut-off function $\mu\in C^\infty(Q)$ such that $0\le\mu\le1$ and
\begin{equation}\label{eq-def-mu}
\mu=\begin{cases}
1 & \mbox{in }Q_{2\de},\\
0 & \mbox{in }Q_0\setminus Q_\de.
\end{cases}
\end{equation}
Then it follows from \eqref{eq-asp-de} that
\begin{equation}\label{eq-van-mu}
\mu(x,0)=1\quad\mbox{for }x\in\ov\Om\,.
\end{equation}
On the other hand, we calculate
\[\pa_t^2(\mu w)=\mu\,\pa_t^2w+[\pa_t^2,\mu]w,\quad\cA(t)(\mu w)=\mu\,\cA(t)w+[\cA(t),\mu]w,\]
where
\[[\pa_t^2,\mu]w:=2(\pa_t\mu)\pa_tw+(\pa_t^2\mu)w,\quad[\cA(t),\mu]w:=2\,a\nb\mu\cdot\nb w+\left\{\rdiv(a\nb\mu)+b\cdot\nb\mu\right\}w.\]
Parallel calculations for $\cA'(t)$ and $\cA''(t)$ yields
\begin{equation}\label{eq-commut1}
\cA'(t)(\mu w)=\mu\,\cA'(t)w+[\cA'(t),\mu]w,\quad\cA''(t)(\mu w)=\mu\,\cA''(t)w+[\cA''(t),\mu]w,
\end{equation}
where
\begin{align*}
[\cA'(t),\mu]w & :=2(\pa_ta)\nb\mu\cdot\nb w+\left\{\rdiv((\pa_ta)\nb\mu)+(\pa_tb)\cdot\nb\mu\right\}w,\\
[\cA''(t),\mu]w & :=2(\pa_t^2a)\nb\mu\cdot\nb w+\left\{\rdiv((\pa_t^2a)\nb\mu)+(\pa_t^2b)\cdot\nb\mu\right\}w.
\end{align*}
Introducing $[\cH,\mu]:=[\pa_t^2,\mu]-[\cA(t),\mu]$, we obtain
\begin{equation}\label{eq-commut0}
\cH(\mu w)=\mu\,\cH w+[\cH,\mu]w.
\end{equation}
Since $[\cH,\mu]$, $[\cA'(t),\mu]$ and $[\cA''(t),\mu]$ are first-order differential operators which only involve derivatives of $\mu$ as coefficients, the definition \eqref{eq-def-mu} of $\mu$ implies
\begin{equation}\label{eq-van}
[\cH,\mu]w=[\cA'(t),\mu]w=[\cA''(t),\mu]w=0\quad\mbox{in }Q_{2\de}\cup(Q_0\setminus Q_\de).
\end{equation}

We set $z:=\mu y$, $z_1:=\mu y_1$ and $z_2:=\mu y_2$. By the definition \eqref{eq-def-mu} of $\mu$, we can take the zero extensions of $z,z_1,z_2$ so that they are defined in $Q$. Therefore, applying formulas \eqref{eq-commut0} and \eqref{eq-commut1} to \eqref{eq-gov-y}, we deduce
\begin{align}
\cH z & =\mu\,\cH y+[\cH,\mu]y=\mu(fR-\cH\wt u\,)+[\cH,\mu]y=:F_0,\label{eq-gov-z}\\
\cH z_1 & =\mu\,\cH y_1+[\cH,\mu]y_1=\mu(\cA'(t)y+f\,\pa_tR-\pa_t(\cH\wt u\,))+[\cH,\mu]y_1\nonumber\\
& =\cA'(t)z+\mu(f\,\pa_tR-\pa_t(\cH\wt u\,))+[\cH,\mu]y_1-[\cA'(t),\mu]y=:F_1,\label{eq-gov-z1}\\
\cH z_2 & =\mu\,\cH y_2+[\cH,\mu]y_2=\mu(2\cA'(t)y_1+\cA''(t)y+f\,\pa_t^2R-\pa_t^2(\cH\wt u\,))+[\cH,\mu]y_2\nonumber\\
& =2\cA'(t)z_1+\cA''(t)z+\mu(f\,\pa_t^2R-\pa_t^2(\cH\wt u\,))\nonumber\\
& \quad\,+[\cH,\mu]y_2-2[\cA'(t),\mu]y_1-[\cA''(t),\mu]y.\label{eq-gov-z2}
\end{align}

Since $y\in\bigcap_{k=0}^2H^{4-k}(-T,T;H^k(\Om))$ and $\mu\in C_0^\infty(Q_0)$ imply $z,\pa_tz\in H^2(Q)$ and $\supp\,z\subset Q_0$, we can apply the Carleman estimate \eqref{eq-CE-1} in Lemma \ref{lem-Carleman} and \eqref{eq-CE-0} to \eqref{eq-gov-z} to obtain
\begin{align}
\sum_{|\ga|\le2}\|(\pa_x^\ga z)\,\e^{s\vp}\|_{L^2(Q_0)}^2 & =\int_{Q_0}\left\{\sum_{i,j=1}^n|\pa_i\pa_jz|^2+|\nb z|^2+|z|^2\right\}\e^{2s\vp}\,\rd x\rd t\nonumber\\
& \le C\int_{Q_0}\left(\f1s|\pa_tF_0|^2+s\,|F_0|^2\right)\e^{2s\vp}\,\rd x\rd t\label{eq-CE-z}
\end{align}
for all large $s>0$. Henceforth we set $\mu_k:=\e^{k\la\de}$ for $k=1,2,3$. Then we employ \eqref{eq-est-wtu}, \eqref{eq-van} and the property of the weight function to estimate
\begin{align}
\int_{Q_0}|F_0|^2\,\e^{2s\vp}\,\rd x\rd t & \le3\int_{Q_0}\left(|\mu fR|^2+|\mu\,\cH\wt u|^2+|[\cH,\mu]y|^2\right)\e^{2s\vp}\,\rd x\rd t\nonumber\\
& \le3\|fR\,\e^{s\vp}\|_{L^2(Q_0)}^2+3\exp\left(2s\max_{\ov Q}\vp\right)\|\cH\wt u\|_{L^2(Q_0)}^2\nonumber\\
& \quad\,+3\int_{Q_\de\setminus Q_{2\de}}|[\cH,\mu]y|^2\,\e^{2s\vp}\,\rd x\rd t\nonumber\\
& \le C\|f\,\e^{s\vp}\|_{L^2(Q_0)}^2+C\,\e^{Cs}\|\wt u\|_{H^2(Q)}^2+3\exp\left(2s\max_{\ov{Q_\de}\setminus Q_{2\de}}\vp\right)\|y\|_{H^1(Q)}^2\nonumber\\
& \le C\|f\,\e^{s\vp}\|_{L^2(Q_0)}^2+C\,\e^{Cs}D^2+C\,\e^{2\mu_2s}M^2\label{eq-est-F0}
\end{align}
for all large $s>0$. Similarly, by
\[\pa_tF_0=(\pa_t\mu)fR+\mu f\,\pa_tR-(\pa_t\mu)\cH\wt u-\mu\,\pa_t(\cH\wt u\,)+\pa_t([\cH,\mu]y)\]
we estimate
\begin{align}
\int_{Q_0}|\pa_tF_0|^2\,\e^{2s\vp}\,\rd x\rd t & \le C\|f\,\e^{s\vp}\|_{L^2(Q_0)}^2+C\,\e^{Cs}\left(\|\wt u\|_{H^3(-T,T;L^2(\Om))}+\|\wt u\|_{H^1(-T,T;H^2(\Om))}\right)\nonumber\\
& \quad\,+C\,\e^{2\mu_2s}\|y\|_{H^2(Q)}^2\nonumber\\
& \le C\|f\,\e^{s\vp}\|_{L^2(Q_0)}^2+C\,\e^{Cs}D^2+C\,\e^{2\mu_2s}M^2\label{eq-est-F0t}
\end{align}
for all large $s>0$. Substituting \eqref{eq-est-F0}--\eqref{eq-est-F0t} into \eqref{eq-CE-z}, we obtain
\begin{equation}\label{eq-est-z}
\sum_{|\ga|\le2}\|(\pa_x^\ga z)\,\e^{s\vp}\|_{L^2(Q_0)}^2\le C\,s\|f\,\e^{s\vp}\|_{L^2(Q_0)}^2+C\,\e^{Cs}D^2+C\,s\,\e^{2\mu_2s}M^2.
\end{equation}

By a similar argument as that for \eqref{eq-gov-z}, we apply the Carleman estimates \eqref{eq-CE-1} and \eqref{eq-CE-0} to \eqref{eq-gov-z1} to estimate
\begin{equation}\label{eq-CE-z1}
\sum_{|\ga|\le2}\|(\pa_x^\ga z_1)\,\e^{s\vp}\|_{L^2(Q_0)}^2\le C\int_{Q_0}\left(\f1s|\pa_tF_1|^2+s\,|F_1|^2\right)\e^{2s\vp}\,\rd x\rd t
\end{equation}
for all large $s>0$, where we utilize \eqref{eq-est-z} to dominate the term with $|F_1|^2$ as
\begin{align}
\int_{Q_0}|F_1|^2\,\e^{2s\vp}\,\rd x\rd t & \le C\int_{Q_0}|\cA'(t)z|^2\,\e^{2s\vp}\,\rd x\rd t+C\|f\,\e^{s\vp}\|_{L^2(Q_0)}^2\nonumber\\
& \quad\,+C\,\e^{Cs}\left(\|\wt u\|_{H^3(-T,T;L^2(\Om))}^2+\|\wt u\|_{H^1(-T,T;H^2(\Om))}^2\right)+C\,\e^{2\mu_2s}\|y\|_{H^2(Q)}^2\nonumber\\
& \le C\sum_{|\ga|\le2}\|(\pa_x^\ga z)\,\e^{s\vp}\|_{L^2(Q_0)}^2+C\|f\,\e^{s\vp}\|_{L^2(Q_0)}^2+C\,\e^{Cs}D^2+C\,\e^{2\mu_2s}M^2\nonumber\\
& \le C\,s\|f\,\e^{s\vp}\|_{L^2(Q_0)}^2+C\,\e^{Cs}D^2+C\,s\,\e^{2\mu_2s}M^2\label{eq-est-F1}
\end{align}
for all large $s>0$. To treat $|\pa_tF_1|^2$, we employ \eqref{eq-commut1} to find
\[\pa_t(\cA'(t)z)=\cA''(t)z+\cA'(t)(\pa_tz)=\cA''(t)z+\cA'(t)z_1+\cA'(t)((\pa_t\mu)y)\]
and thus
\begin{align}
& \quad\,\,\int_{Q_0}|\pa_tF_1|^2\,\e^{2s\vp}\,\rd x\rd t\nonumber\\
& \le C\sum_{|\ga|\le2}\left(\|(\pa_x^\ga z_1)\,\e^{s\vp}\|_{L^2(Q_0)}^2+\|(\pa_x^\ga z)\,\e^{s\vp}\|_{L^2(Q_0)}^2\right)+C\|f\,\e^{s\vp}\|_{L^2(Q_0)}^2\nonumber\\
& \quad\,+C\,\e^{Cs}\left(\|\wt u\|_{H^4(-T,T;L^2(\Om))}^2+\|\wt u\|_{H^2(-T,T;H^2(\Om))}^2\right)+C\,\e^{2\mu_2s}\sum_{k=0}^2\|y\|_{H^{3-k}(-T,T;H^k(\Om))}^2\nonumber\\
& \le C\sum_{|\ga|\le2}\|(\pa_x^\ga z_1)\,\e^{s\vp}\|_{L^2(Q_0)}^2+C\,s\|f\,\e^{s\vp}\|_{L^2(Q_0)}^2+C\,\e^{Cs}D^2+C\,s\,\e^{2\mu_2s}M^2\label{eq-est-F1t}
\end{align}
for all large $s>0$. Substituting \eqref{eq-est-F1}--\eqref{eq-est-F1t} into \eqref{eq-CE-z1} yields
\begin{align*}
\sum_{|\ga|\le2}\|(\pa_x^\ga z_1)\,\e^{s\vp}\|_{L^2(Q_0)}^2 & \le\f Cs\sum_{|\ga|\le2}\|(\pa_x^\ga z_1)\,\e^{s\vp}\|_{L^2(Q_0)}^2+C\,s^2\|f\,\e^{s\vp}\|_{L^2(Q_0)}^2\\
& \quad\,+C\,\e^{Cs}D^2+C\,s^2\,\e^{2\mu_2s}M^2.
\end{align*}
for all large $s>0$. Choosing $s>0$ large enough and absorbing the term $\f1s\sum_{i,j=1}^n|\pa_i\pa_jz_1|^2$ on the right-hand side into the left-hand side, we have for all large $s>0$ that
\begin{equation}\label{eq-est-z1}
\sum_{|\ga|\le2}\|(\pa_x^\ga z_1)\,\e^{s\vp}\|_{L^2(Q_0)}^2\le C\,s^2\|f\,\e^{s\vp}\|_{L^2(Q_0)}^2+C\,\e^{Cs}D^2+C\,s^2\,\e^{2\mu_2s}M^2.
\end{equation}

Finally, applying the first-order Carleman estimate \eqref{eq-CE-0} to \eqref{eq-gov-z2}, we use estimates \eqref{eq-est-z} and \eqref{eq-est-z1} to deduce
\begin{align}
& \quad\,\,s\int_{Q_0}\left(|\pa_tz_2|^2+s^2|z_2|^2\right)\e^{2s\vp}\,\rd x\rd t\nonumber\\
& \le C\sum_{|\ga|\le2}\left(\|(\pa_x^\ga z_1)\,\e^{s\vp}\|_{L^2(Q_0)}^2+\|(\pa_x^\ga z)\,\e^{s\vp}\|_{L^2(Q_0)}^2\right)+C\|f\,\e^{s\vp}\|_{L^2(Q_0)}^2\nonumber\\
& \quad\,+C\,\e^{Cs}\left(\|\wt u\|_{H^4(-T,T;L^2(\Om))}^2+\|\wt u\|_{H^2(-T,T;H^2(\Om))}^2\right)+C\,\e^{2\mu_2s}\sum_{k=0}^1\|y\|_{H^{3-k}(-T,T;H^k(\Om))}^2\nonumber\\
& \le C\,s^2\|f\,\e^{s\vp}\|_{L^2(Q_0)}^2+C\,\e^{Cs}D^2+C\,s^2\,\e^{2\mu_2s}M^2.\label{eq-CE-z2}
\end{align}
for all large $s>0$. Meanwhile, by $\supp\,z_2\subset Q_0$, we obtain
\begin{align}
\|(z_2\,\e^{s\vp})(\,\cdot\,,0)\|_{L^2(\Om)}^2 & =\int_{-T}^0\left(\f\rd{\rd t}\|(z_2\,\e^{s\vp})(\,\cdot\,,t)\|_{L^2(\Om)}^2\right)\rd t\nonumber\\
& =2\int_{-T}^0\!\int_\Om\left\{(\pa_tz_2)z_2+s(\pa_t\vp)|z_2|^2\right\}\e^{2s\vp}\,\rd x\rd t\nonumber\\
& \le C\int_{Q_0}\left(|\pa_tz_2||z_2|+|z_2|^2\right)\e^{2s\vp}\,\rd x\rd t\nonumber\\
& \le C\int_{Q_0}\left(\f1s|\pa_tz_2|^2+s\,|z_2|^2\right)\e^{2s\vp}\,\rd x\rd t\label{eq-ineq1}\\
& \le C\|f\,\e^{s\vp}\|_{L^2(Q_0)}^2+C\,\e^{Cs}D^2+C\,\e^{2\mu_2s}M^2\label{eq-ineq2}
\end{align}
for all large $s>0$, where \eqref{eq-ineq1} and \eqref{eq-ineq2} follow from the inequality
\[|\pa_tz_2||z_2|=\f1{\sqrt s}|\pa_tz_2|\cdot\sqrt s\,|z_2|\le\f1{2s}|\pa_tz_2|^2+\f s2|z_2|^2\]
and the estimate \eqref{eq-CE-z2}, respectively. On the other hand, taking $t=0$ in the governing equation of $y$ in \eqref{eq-gov-y}, we use \eqref{eq-van-mu} and \eqref{eq-IC-u} to find
\begin{align*}
z_2(\,\cdot\,,0) & =(\mu y_2)(\,\cdot\,,0)=y_2(\,\cdot\,,0)=\pa_t^2y(\,\cdot\,,0)=(\cA(t)(u-\wt u)+fR-\cH\wt u\,)(\,\cdot\,,0)\\
& =(fR-\pa_t^2\wt u\,)(\,\cdot\,,0).
\end{align*}
The above equality holds in the sense of $L^2(\Om)$ because $\pa_t^2y\in H^2(-T,T;L^2(\Om))$ and $R\in W^{2,\infty}(-T,T;L^\infty(\Om))$. Recalling the assumption \eqref{eq-asp-R}, we apply \eqref{eq-ineq2} to further estimate for all large $s>0$ that
\begin{align*}
\|(f\,\e^{s\vp})(\,\cdot\,,0)\|_{L^2(\Om)}^2 & \le C\|(fR\,\e^{s\vp})(\,\cdot\,,0)\|_{L^2(\Om)}^2\\
& \le C\|(z_2\,\e^{s\vp})(\,\cdot\,,0)\|_{L^2(\Om)}^2+C\|((\pa_t^2\wt u\,)\,\e^{s\vp})(\,\cdot\,,0)\|_{L^2(\Om)}^2\\
& \le C\|f\,\e^{s\vp}\|_{L^2(Q_0)}^2+C\,\e^{Cs}D^2+C\,\e^{2\mu_2s}M^2+C\,\e^{Cs}\|\wt u\|_{H^3(-T,T;L^2(\Om))}^2\\
& \le C\|f\,\e^{s\vp}\|_{L^2(Q_0)}^2+C\,\e^{Cs}D^2+C\,\e^{2\mu_2s}M^2.
\end{align*}
Now we treat the term $\|f\,\e^{s\vp}\|_{L^2(Q_0)}^2$ as
\begin{align*}
\|f\,\e^{s\vp}\|_{L^2(Q_0)}^2 & =\int_\Om|f|^2\left(\int^T_{-T}\e^{2s\vp(x,t)}\,\rd t\right)\rd x\\
& =\int_\Om|f|^2\,\e^{2s\vp(x,0)}\left(\int^T_{-T}\e^{-2s(\vp(x,0)-\vp(x,t))}\,\rd t\right)\rd x.
\end{align*}
Since $\vp(x,0)>\vp(x,t)$ for $t\ne0$, Lebesgue's dominated convergence theorem yields
\[\int^T_{-T}\e^{-2s(\vp(x,0)-\vp(x,t))}\,\rd t=o(1)\quad\mbox{as }s\to\infty.\]
Since $C>0$ is independent of $s$, this indicates
\[\|(f\,\e^{s\vp})(\,\cdot\,,0)\|_{L^2(\Om)}^2\le o(1)\|(f\,\e^{s\vp})(\,\cdot\,,0)\|_{L^2(\Om)}^2+C\,\e^{Cs}D^2+C\,\e^{2s\mu_2}M^2\]
for all large $s\ge s_*$, where $s_*>0$ is a large constant. We can choose $s_*>0$ large again, so that we can absorb the first term on the right-hand side into the left-hand side to obtain
\begin{equation}\label{eq-est-f1}
\|(f\,\e^{s\vp})(\,\cdot\,,0)\|_{L^2(\Om)}^2\le C\,\e^{Cs}D^2+C\,\e^{2s\mu_2}M^2
\end{equation}
for all $s\ge s_*$. On the other hand, we estimate $\|(f\,\e^{s\vp})(\,\cdot\,,0)\|_{L^2(\Om)}^2$ by
\begin{equation}\label{eq-est-f2}
\|(f\,\e^{s\vp})(\,\cdot\,,0)\|_{L^2(\Om)}^2\ge\|(f\,\e^{s\vp})(\,\cdot\,,0)\|_{L^2(\Om_{3\de})}^2\ge\e^{2\mu_3s}\|f\|_{L^2(\Om_{3\de})}^2,
\end{equation}
where $\mu_3=\e^{3\la\de}$ and $\Om_{3\de}$ was defined as that in \eqref{eq-def-Qde}. Inequalities \eqref{eq-est-f1} and \eqref{eq-est-f2} yield
\begin{equation}\label{eq-est-f3}
\|f\|_{L^2(\Om_{3\de})}^2\le C\,\e^{Cs}D^2+C\,\e^{-2s(\mu_3-\mu_2)}M^2
\end{equation}
for all $s\ge s_*$. Replacing $C$ by $C\,\e^{Cs_*}$, we have \eqref{eq-est-f3} for all $s\ge0$. We note that $\mu_3-\mu_2=\e^{3\la\de}-\e^{2\la\de}>0$.

We consider the two cases $D\ge M$ and $D<M$ separately.

{\bf Case 1 } If $D\ge M$, then \eqref{eq-est-f3} directly implies
\begin{equation}\label{eq-est-f4}
\|f\|^2_{L^2(\Om_{3\de})}\le 2C\,\e^{Cs}D^2.
\end{equation}

{\bf Case 2 } If $D<M$, we can choose $s>0$ suitably to minimize the right-hand side of \eqref{eq-est-f3} such that
\[\e^{Cs}D^2=\e^{-2s(\mu_3-\mu_2)}M^2,\]
that is,
\[s=\f2{C+2(\mu_3-\mu_2)}\log\f MD>0.\]
Then we obtain
\begin{equation}\label{eq-est-f5}
\|f\|^2_{L^2(\Om_{3\de})}\le 2CM^{2(1-\te)}D^{2\te},
\end{equation}
where
\[\te=\f{2(\mu_3-\mu_2)}{C+2(\mu_3-\mu_2)}\in(0,1).\]

Finally, replacing $3\de$ by $\de$, we see that estimates \eqref{eq-est-f4} and \eqref{eq-est-f5} yield \eqref{eq-thm-est}. Thus the proof of Theorem \ref{thm-ISP} is completed.
\end{proof}

\Section{Iterative Thresholding Algorithm}\label{sec-recon}

Based on the theoretical stability explained in the previous sections, this section aims at the development of an effective algorithm for the numerical reconstruction of the source term.

For conciseness, we consider the initial-boundary value problem for a wave equation with the homogeneous Neumann boundary condition
\begin{equation}\label{eq-gov-w}
\begin{cases}
\pa_t^2u(x,t)=\tri u(x,t)+f(x)R(x,t) & (x\in\Om,\ 0<t<T),\\
u(x,0)=\pa_tu(x,0)=0 & (x\in\Om),\\
\pa_\nu u(x,t)=0 & (x\in\pa\Om,\ 0<t<T),
\end{cases}
\end{equation}
whose solution will be denoted as $u(f)$ to emphasize its dependency upon the unknown component $f$. As a typical situation, we only treat Case (II) of Problem \ref{prob-ISP}, namely, the determination of $f$ by the partial interior observation data $u(f)|_{\om\times(0,T)}$ with a subdomain $\om\subset\Om$. Except for its simplicity, we restrict our discussion to the time-independent formulation \eqref{eq-gov-w} instead of the time-dependent one because the underlying ill-posedness are in principle the same. Moreover, in this case a global Lipschitz stability is guaranteed by Lemma \ref{lem-inter}, and a necessary condition on the geometry of observable regions is given explicitly in \eqref{eq-cond}.

For later use, first we give a definition of the generalized solution to \eqref{eq-gov-w}.

\begin{defi}[see Isakov \cite{I06}]\label{def-sol}
Let $f\in L^2(\Om)$ and $R\in L^2(0,T;L^\infty(\Om))$. We say that $u(f)\in H^1(\Om\times(0,T))$ is a generalized solution to problem \eqref{eq-gov-w} if it satisfies
\[\int_0^T\!\!\!\int_\Om(\nb u(f)\cdot\nb v-(\pa_tu(f))\pa_tv)\,\rd x\rd t=\int_0^T\!\!\!\int_\Om fR\,v\,\rd x\rd t\]
for any test function $v\in H^1(\Om\times(0,T))$ with $v|_{t=T}=0$, and the initial condition $u(f)|_{t=0}=0$.
\end{defi}

The above definition of the generalized solution is easily understood by performing integration by parts to sufficiently smooth solutions. Furthermore, it is in accordance with the classical well-posedness result (see Lemma \ref{lem-Lions}(a) with $k=1$).

Henceforth, we specify $f_\true\in L^2(\Om)$ as the true solution to Problem \ref{prob-ISP} and suppose that we are provided the noise contaminated observation data $u^\de$ in $\om\times(0,T)$ satisfying $\|u^\de-u(f_\true)\|_{L^2(\om\times(0,T))}\le\de$, where $\de>0$ stands for the noise level. For avoiding ambiguity, we interpret $u^\de\equiv0$ out of $\om\times(0,T)$ so that it is well-defined in $\Om\times(0,T)$.

With the a priori knowledge on the boundedness of $f_\true$ and appropriate observation data, the reconstruction can be carried out through a classical Tikhonov regularization technique. We formulate the reconstruction as the following output least squares formulation with the Tikhonov regularization
\begin{equation}\label{eq-n3}
\min_{f\in L^2(\Om)}J(f),\quad J(f):=\|u(f)-u^\de\|_{L^2(\om\times(0,T))}^2+\al\|f\|_{L^2(\Om)}^2,
\end{equation}
where $\al>0$ is the regularization parameter.

Nearly all effective iterative methods for solving nonlinear optimizations need the information of the derivatives of the concerned objective functional. It follows from a direct computation that the Fr\'echet derivative $J'(f)g$ of $J(f)$ for any direction $g\in L^2(\Om)$ reads
\begin{align}
J'(f)g & =2\int_0^T\!\!\!\int_\om\left(u(f)-u^\de\right)(u'(f)g)\,\rd x\rd t+2\al\int_\Om fg\,\rd x\nonumber\\
& =2\int_0^T\!\!\!\int_\om\left(u(f)-u^\de\right)u(g)\,\rd x\rd t+2\al\int_\Om fg\,\rd x.\label{eq-nn1}
\end{align}
Here $u'(f)g$ denotes the Fr\'echet derivative of $u(f)$ in the direction $g$, and the linearity of \eqref{eq-gov-w} immediately yields
\[u'(f)g=\lim_{\ep\to0}\f{u(f+\ep\,g)-u(f)}\ep=u(g).\]
Obviously, it is extremely expensive to use this formula to evaluate $J'(f)g$ for all $g\in L^2(\Om)$, since one should solve system \eqref{eq-gov-w} for $u(g)$ with $g$ varying in $L^2(\Om)$ in the computation for a fixed $f$.

In order to reduce the computational costs for computing the Fr\'echet derivatives, we introduce the adjoint system of \eqref{eq-gov-w}, that is, the following system for a backward wave equation
\begin{equation}\label{eq-gov-v}
\begin{cases}
\pa_t^2v-\tri v=\chi_\om\left(u(f)-u^\de\right) & \mbox{in }\Om\times(0,T),\\
v=\pa_tv=0 & \mbox{in }\Om\times\{0\},\\
\pa_\nu v=0 & \mbox{on }\pa\Om\times(0,T).
\end{cases}
\end{equation}
Here $\chi_\om$ is the characterization function of $\om$, and we shall denote the solution to \eqref{eq-gov-v} as $v(f)$. The generalized solution to \eqref{eq-gov-v} can be defined in the same way as that in Definition \ref{def-sol}. On the other hand, since Lemma \ref{lem-Lions} and $f\in L^2(\Om)$ indicate
\begin{equation}\label{eq-def-D}
u(f)\in\cD:=\left\{w\in C([0,T];H^1(\Om));\,\pa_tw\in C([0,T];L^2(\Om))\right\},
\end{equation}
we have $\chi_\om\left(u(f)-u^\de\right)\in L^2(\Om\times(0,T))$ and thus $v(f)\in\cD$ again by Lemma \ref{lem-Lions}. Therefore, for any $f,g\in L^2(\Om)$, regarding $u(g)$ and $v(f)$ as the test functions of each other, we can further treat the first term in \eqref{eq-nn1} as
\begin{align}
\int_0^T\!\!\!\int_\om\left(u(f)-u^\de\right)u(g)\,\rd x\rd t & =\int_0^T\!\!\!\int_\Om\chi_\om\left(u(f)-u^\de\right)u(g)\,\rd x\rd t\nonumber\\
& =\int_0^T\!\!\!\int_\Om\left(\nb u(g)\cdot\nb v(f)-(\pa_tu(g))\pa_tv(f)\right)\rd x\rd t\nonumber\\
& =\int_0^T\!\!\!\int_\Om gR\,v(f)\,\rd x\rd t,\label{eq-n5}
\end{align}
implying
\[J'(f)g=2\int_\Om\left(\int_0^TR\,v(f)\,\rd t+\al\,f\right)g\,\rd x,\quad g\in L^2(\Om).\]
This suggests a characterization of the solution to the minimization problem \eqref{eq-n3}.

\begin{prop}\label{lem-min}
$f^*\in L^2(\Om)$ is the minimizer of the functional $J(f)$ in \eqref{eq-n3} only if it satisfies the Euler equation
\begin{equation}\label{eq-min-f}
\int_0^TR\,v(f^*)\,\rd t+\al\,f^*=0,
\end{equation}
where $v(f^*)$ solves the backward system \eqref{eq-gov-v} with the coefficient $f^*$.
\end{prop}

To solve the nonlinear equation \eqref{eq-min-f} for $f^*$, one may employ the iteration
\begin{equation}\label{eq-itr-f}
f_{m+1}=\f K{K+\al}f_m-\f1{K+\al}\int_0^TR\,v(f_m)\,\rd t,
\end{equation}
where $K>0$ is a tuning parameter acting as a weight between the previous step and the iterative update.

To discuss the choice of $K$ to guarantee the convergence, we take advantage of the fact that the iteration \eqref{eq-itr-f} in principle coincides with the iterative thresholding algorithm, which can be derived from the minimization problem of a surrogate functional (see, e.g. \cite{DDM04}). Actually, fixing $g\in L^2(\Om)$, we introduce a surrogate functional $J^s(f,g)$ of $J(f)$ as
\[J^s(f,g):=J(f)+K\|f-g\|_{L^2(\Om)}^2-\|u(f)-u(g)\|_{L^2(\om\times(0,T))}^2.\]
For the positivity of $J^s$, there should hold $K\|f\|_{L^2(\Om)}^2\ge\|u(f)\|_{L^2(\om\times(0,T))}^2$ for all $f\in L^2(\Om)$. This is achieved by choosing
\begin{equation}\label{eq-coef-M}
K\ge\|A\|^2,
\end{equation}
where $A$ is a linear operator defined as
\begin{align*}
A:L^2(\Om) & \to L^2(\om\times(0,T)),\\
f & \mapsto u(f)|_{\om\times(0,T)},
\end{align*}
and the boundedness of $A$ is readily seen from Lemma \ref{lem-Lions}. Therefore, there holds
\[J(f)=J^s(f,f)\le J^s(f,g),\]
and thus $J^s(f,g)$ can be regarded as a small perturbation of $J(f)$ when $g$ is close to $f$. On the other hand, it follows from \eqref{eq-n5} that
\begin{align*}
J^s(f,g) & =2\int_0^T\!\!\!\int_\om u(f)\left(u(g)-u^\de\right)\rd x\rd t+\|u^\de\|_{L^2(\om\times(0,T))}^2-\|u(g)\|_{L^2(\om\times(0,T))}^2\\
& \quad\,+\al\|f\|_{L^2(\Om)}^2+K\|f-g\|_{L^2(\Om)}^2\\
& =K\|f-g\|_{L^2(\Om)}^2+\al\|f\|_{L^2(\Om)}^2+2\int_0^T\!\!\!\int_\Om fR\,v(g)\,\rd x\rd t\\
& \quad\,-\|u(g)\|_{L^2(\om\times(0,T))}^2+\|u^\de\|_{L^2(\om\times(0,T))}^2\\
& =(K+\al)\|f\|_{L^2(\Om)}^2-2\int_\Om f\left(K\,g-\int_0^TR\,v(g)\,\rd t\right)\rd x\\
& \quad\,+K\|g\|_{L^2(\Om)}^2-\|u(g)\|_{L^2(\om\times(0,T))}^2+\|u^\de\|_{L^2(\om\times(0,T))}^2.
\end{align*}
Since this is a quadratic form with respect to $f$ when $u^\de$ and $g$ are fixed, we see
\[\mathop{\mathrm{arg}\min}_{f\in L^2(\Om)}J^s(f,g)=\f K{K+\al}g-\f1{K+\al}\int_0^TR\,v(g)\,\rd t.\]
Consequently, the iterative update \eqref{eq-itr-f} is equivalent to solving the minimization problem $\min_{f\in L^2(\Om)}J^s(f,g)$ with $g=f_m$. Moreover, the convergence of this iteration was proved in \cite{DDM04} for any bounded linear operator $A$, provided that the constant $K>0$ is chosen according to condition \eqref{eq-coef-M}.

Now we are well prepared to state the main algorithm for the numerical reconstruction.

\begin{algo}\label{al-soga}
Choose a tolerance $\ve>0$, a regularization parameter $\al>0$ and a tuning constant $K>0$ according to \eqref{eq-coef-M}. Give an initial guess $f_0$, and set $m=0$.
\begin{enumerate}
\item Compute $f_{m+1}$ by the iterative update \eqref{eq-itr-f}.
\item If $\|f_{m+1}-f_m\|_{L^2(\Om)}/\|f_m\|_{L^2(\Om)}\le\ve$, then stop the iteration. Otherwise, update $m\leftarrow m+1$ and return to Step 1.
\end{enumerate}
\end{algo}

It turns out that at each iteration step, we only need to solve the forward system \eqref{eq-gov-w} once for $u(f_m)$ and the backward system \eqref{eq-gov-v} once for $v(f_m)$ subsequently. Therefore, it is very easy and cheap to implement Algorithm \ref{al-soga}. As will be shown from many numerical experiments in the next section, we see that our proposed Algorithm \ref{al-soga} is also considerably efficient and accurate even for three spatial dimensions.

We conclude this section by stating the convergence result of Algorithm \ref{al-soga}, which is a direct application of \cite[Theorem 3.1]{DDM04}.

\begin{lem}\label{lem-conv}
Let $K>0$ be a constant satisfying condition \eqref{eq-coef-M}. Then for any $f_0\in L^2(\Om)$, the sequence $\{f_m\}_{m=1}^\infty$ produced by the iteration \eqref{eq-itr-f} converges strongly to the solution to the minimization problem \eqref{eq-n3}.
\end{lem}

\Section{Numerical Experiments}\label{sec-numer}

In this section, we will apply the established Algorithm \ref{al-soga} to the numerical identification of the spatial component $f$ of the source term in system \eqref{eq-gov-w}. The general settings of the numerical reconstruction are assigned as follows. For simplicity, we take $\Om=(0,1)^n$ ($n=1,2,3$). The duration $T$ may change with respect to the choices of $\om$ according to condition \eqref{eq-cond} which guarantees a reasonable reconstruction. Actually, since an observable subdomain $\om$ has certain thickness in practice, the condition on $T$ can be relaxed to $T>\diam(\Om\setminus\ov\om)$ in numerical tests. Although in Section \ref{sec-recon} the difference between the noiseless data $u(f_\true)$ and the actually observed data $u^\de$ was evaluated in the $L^2(\om\times(0,T))$-norm, here for simplicity we produce $u^\de$ by adding uniform random noises to the noiseless data, i.e.
\[u^\de(x,t)=u(f_\true)(x,t)+\de\,\mathrm{rand}(-1,1),\quad\forall\,x\in\om,\ \forall\,t\in(0,T),\]
where $\mathrm{rand}(-1,1)$ denotes the uniformly distributed random number in $[-1,1]$ and $\de\ge0$ is the noise level. Here we choose $\de$ as a certain portion of the amplitude of the true solution, i.e.
\[\de:=\de_0\max_{\ov\Om\times[0,T]}|u(f_\true)|,\quad0<\de_0<1.\]
As for the various parameters involved in Algorithm \ref{al-soga}, we take $\ve=1\%\,\de_0$ as the tolerance, and choose the regularization parameter as $\al=0.1\%\,\de$. The initial guess $f_0$ is always taken as a constant, which is usually rather inaccurate in the test problems. Finally, the tuning parameter $K>0$ will be chosen according to the size of the subdomain $\om$, the duration $T$, and the known component $R(x,t)$ of the source term.

At each step of the iteration \eqref{eq-itr-f} in all of the numerical experiments, the forward system \eqref{eq-gov-w} and the backward system \eqref{eq-gov-v} are solved by some absolutely stable schemes of the finite difference method. In our implementations, we apply the von Neumann scheme for one-dimensional case and the alternating direction implicit (ADI) method for two- and three-dimensional cases (see \cite{L62,FM65}). It turns out that the ADI method performs efficiently even in three-dimensional case. In fact, it only takes about 5 seconds for a problem of $50^3\times100$ scale. On the other hand, the involved integrals in time are simply approximated by the composite trapezoidal rule.

In what follows, we shall demonstrate the reconstruction method by abundant test examples in one, two and three spatial dimensions. Other than the illustrative figures, we mainly evaluate the numerical performance by the number $M$ of iterations, the relative $L^2$ error
\[\mathrm{err}:=\f{\|f_M-f_\true\|_{L^2(\Om)}}{\|f_\true\|_{L^2(\Om)}}\]
and the elapsed time, where $f_M$ is recognized as the result of the numerical reconstruction.

\Subsection{One-dimensional examples}

In case of $n=1$, we always take $T=1$ and divide the space-time region $\ov\Om\times[0,T]=[0,1]\times[0,1]$ as a $101\times101$ equidistant mesh and test the performance of Algorithm \ref{al-soga} from various aspects. The choice $T=1$ is sufficient because there always holds $\diam(\Om\setminus\ov\om)<1$ whatever $\om$ we set.

\begin{figure}[htbp]\centering
\includegraphics[trim=4mm 4mm 13mm 1mm,clip=true,width=.45\textwidth]{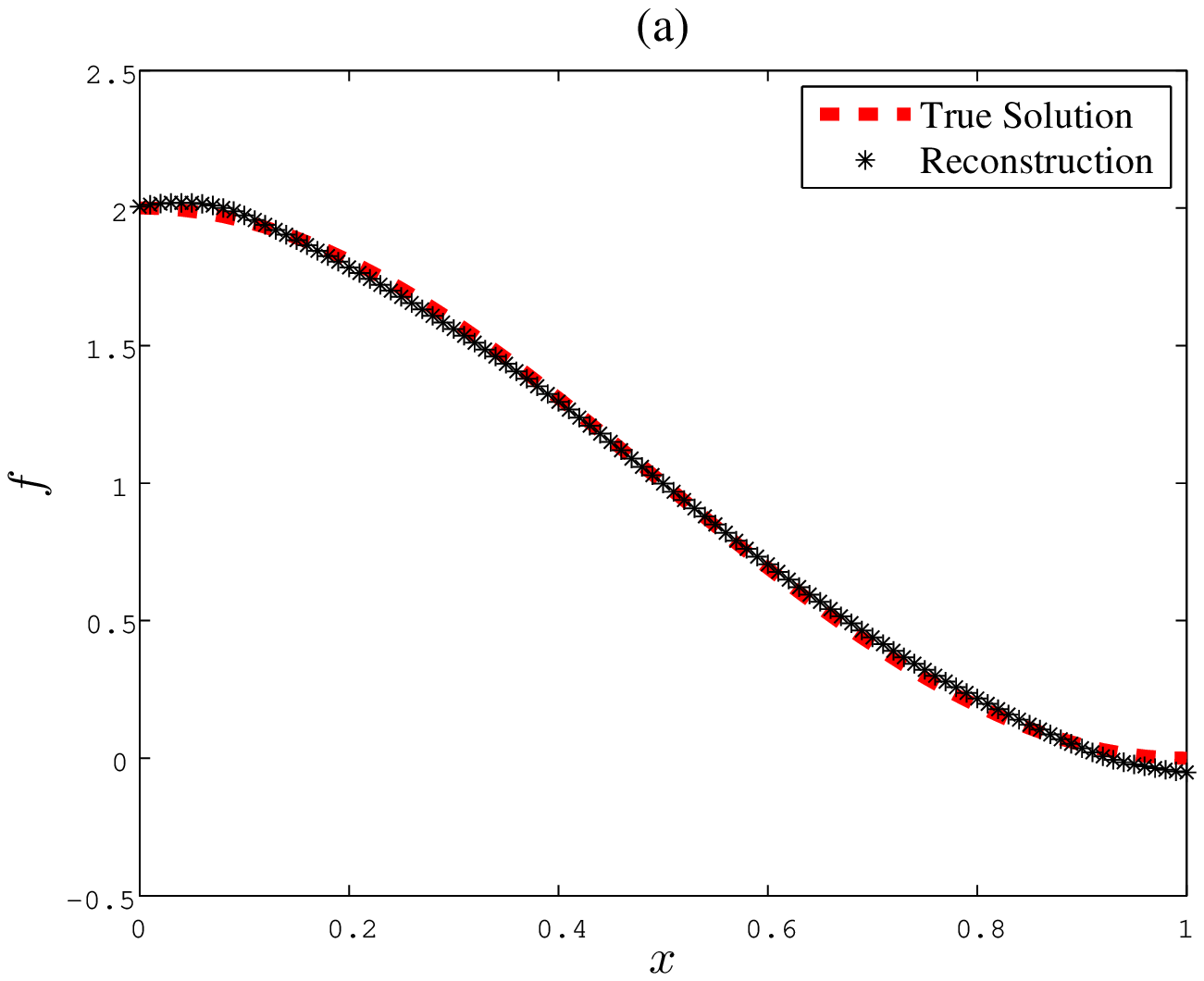}\qquad\includegraphics[trim=4mm 4mm 13mm 1mm,clip=true,width=.45\textwidth]{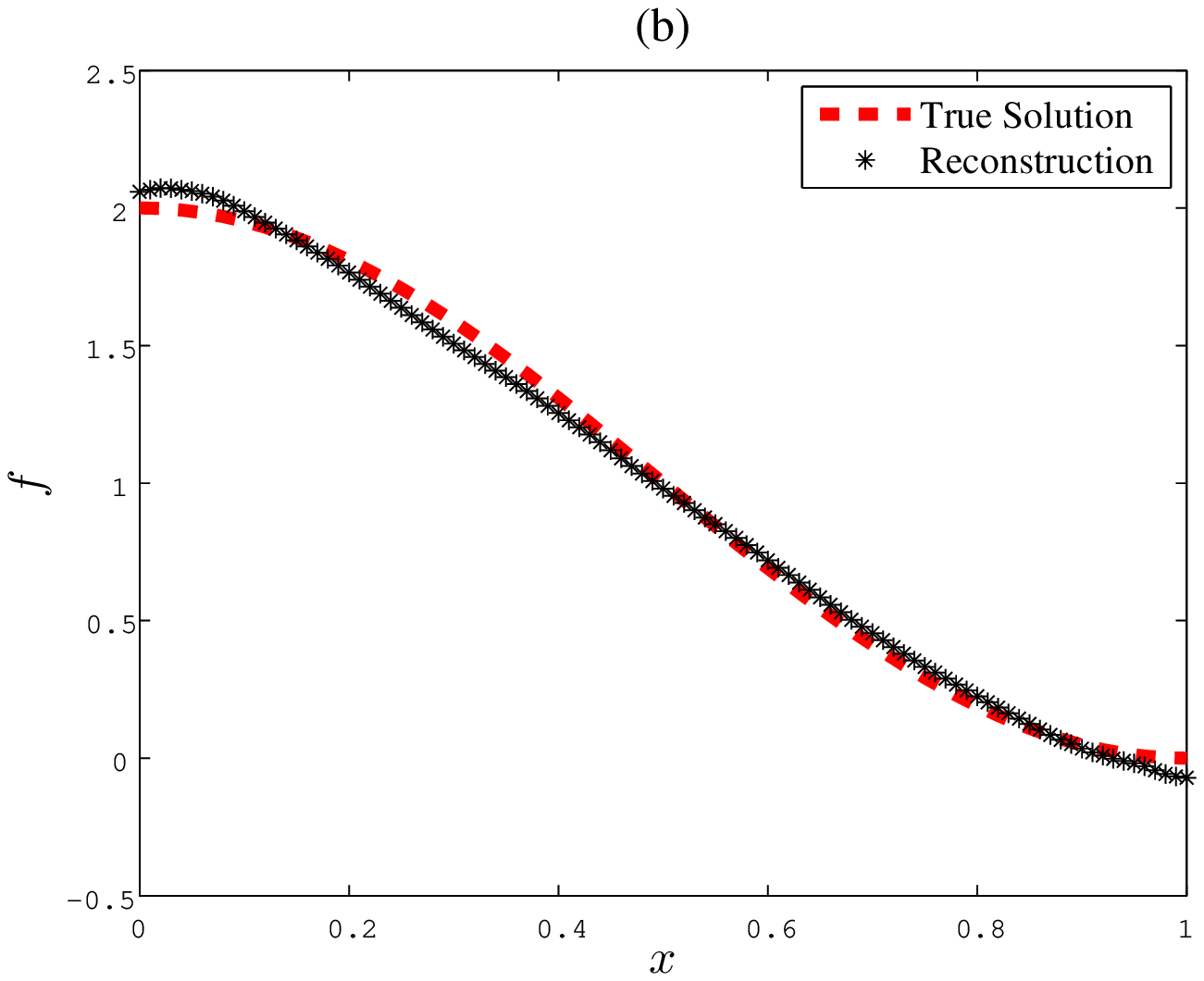}\\
\caption{Illustrations of several reconstructions of $f_\true$ in Example \ref{ex-d1e1} with different choices of the noise level $\de_0$. (a) $\de_0=1\%$. (b) $\de_0=8\%$.}\label{fig-d1e1}
\end{figure}

\begin{ex}\label{ex-d1e1}
In this example, we carry out numerical tests with different combinations of the noise level $\de$ and the observable subdomain $\om$ to see their influences upon the reconstructions. Take the known component of the source term as $R(x,t)=x+t+1$, let $f_\true(x)=\cos(\pi x)+1$ and set the initial guess as $f_0\equiv1$. First we fix $\om=\Om\setminus[0.1,0.9]$ and change the noise levels as $1\%$, $2\%$, $4\%$ and $8\%$ of the amplitude of $u(f_\true)$. Then we fix an $1\%$ noise and reduce the size of $\om$ from $\om=\Om\setminus[0.2,0.8]$ to $\om=\Om\setminus[0.05,0.95]$. The choices of parameters in the tests and corresponding numerical performances are listed in Table \ref{tab-d1e1}. For a better understanding of reconstructions, we visualize several representative examples in Table \ref{tab-d1e1} to compare the true solutions and the recovered ones in Figure \ref{fig-d1e1}.
\end{ex}
\begin{table}[htbp]\centering
\caption{Parameters and corresponding numerical performances in Example \ref{ex-d1e1} under various combinations of noise levels and the observable subdomains.}\label{tab-d1e1}
\begin{tabular}{ccc|ccc|c}
\hline\hline
$\de_0$ & $\om$ & $K$ & $M$ & $\mathrm{err}$ & elapsed time (s) & illustration\\
\hline
$1\%$ & $\Om\setminus[0.1,0.9]$ & $0.02$ & $113$ & $1.86\%$ & $2.71$ & Figure \ref{fig-d1e1}(a)\\
$2\%$ & $\Om\setminus[0.1,0.9]$ & $0.02$ & $84$ & $2.91\%$ & $2.04$\\
$4\%$ & $\Om\setminus[0.1,0.9]$ & $0.02$ & $73$ & $3.32\%$ & $1.65$\\
$8\%$ & $\Om\setminus[0.1,0.9]$ & $0.02$ & $65$ & $3.79\%$ & $1.63$ & Figure \ref{fig-d1e1}(b)\\
\hline
$1\%$ & $\Om\setminus[0.2,0.8]$ & $0.04$ & $118$ & $1.15\%$ & $2.78$\\
$1\%$ & $\Om\setminus[0.05,0.95]$ & $0.015$ & $122$ & $2.77\%$ & $2.81$\\
\hline\hline
\end{tabular}
\end{table}

\begin{ex}\label{ex-d1e2}
Now we compare the numerical performances by selecting various true solutions $f_\true$ with different monotonicity and smoothness. More precisely, we fix $R(x,t)=2+\pi^2t^2$ and choose (a) $f_\true^{\,\a}(x)=x$, (b) $f_\true^{\,\rb}(x)=\sin(\pi x)+x$, (c) $f_\true^{\,\rc}(x)=\f12\cos(2\pi x)+1$ or (d) $f_\true^{\,\rd}(x)=1-|2x-1|$. In all cases, we set the noise level as $5\%$ of the amplitude of $u(f_\true)$, and take $\om=\Om\setminus[0.1,0.9]$. Correspondingly, the tuning parameter is chosen as $K=0.1$. The numbers $M$ of iterations and relative errors are listed in Table \ref{tab-d1e2}. The comparisons of several pairs of the true solutions and the reconstructed ones are shown in Figure \ref{fig-d1e2}.
\begin{table}[htbp]\centering
\caption{Numerical performances of the reconstructions in Example \ref{ex-d1e2} for various choices of true solutions with different smoothness.}\label{tab-d1e2}
\begin{tabular}{cc|cc|c}
\hline\hline
$f_\true(x)$ & initial guess & $M$ & $\mathrm{err}$ & illustration\\
\hline
$x$ & $0.5$ & $6$ & $1.41\%$ & Figure \ref{fig-d1e2}(a)\\
$\sin(\pi x)+x$ & $2.5$ & $43$ & $2.03\%$\\
$\f12\cos(2\pi x)+1$ & $1$ & $179$ & $7.55\%$\\
$1-|2x-1|$ & $0.5$ & $223$ & $11.41\%$ & Figure \ref{fig-d1e2}(d)\\
\hline\hline
\end{tabular}
\end{table}

\begin{figure}[htbp]\centering
\includegraphics[trim=7mm 4mm 13mm 1mm,clip=true,width=.45\textwidth]{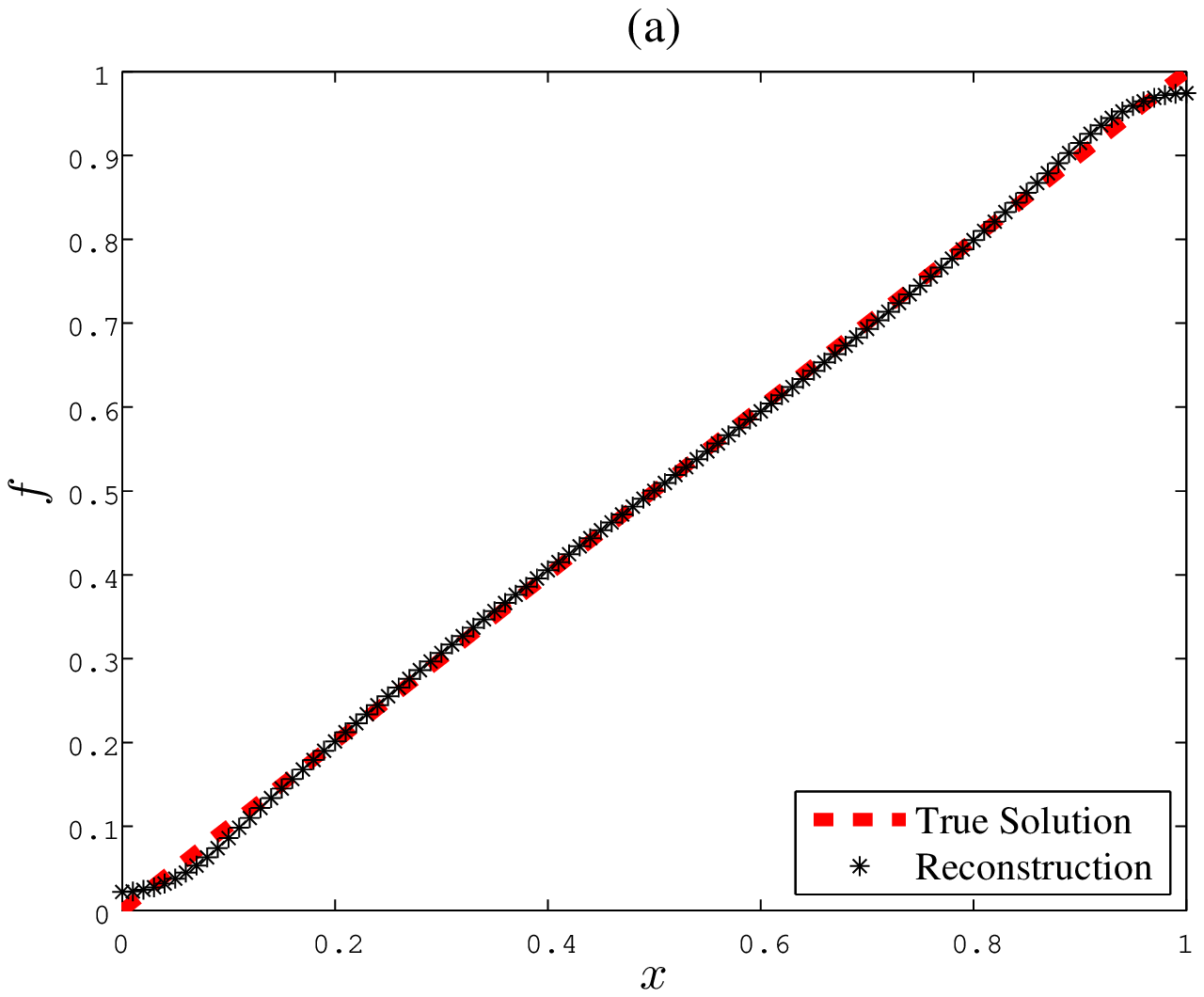}\qquad\includegraphics[trim=7mm 4mm 13mm 1mm,clip=true,width=.45\textwidth]{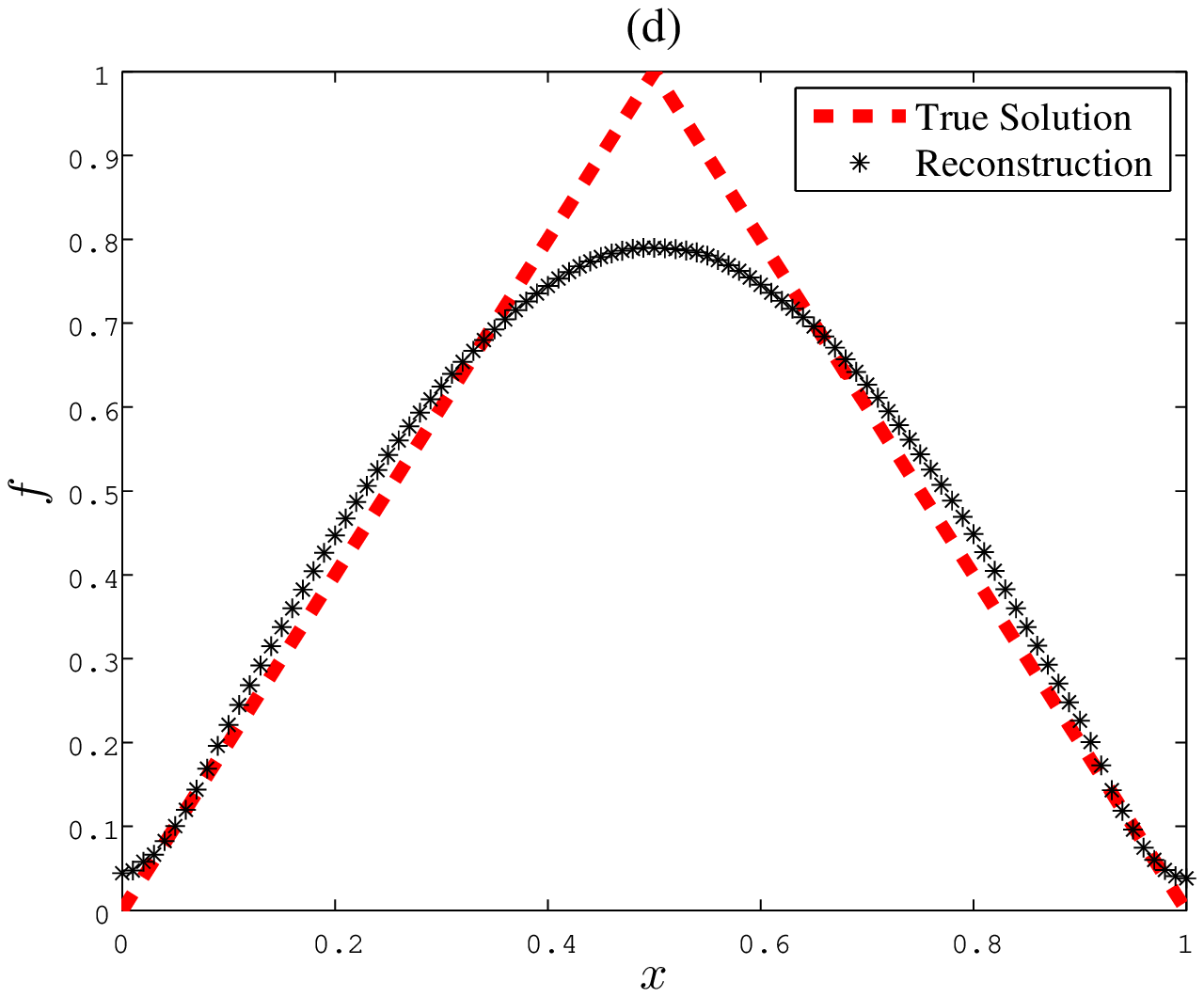}\\
\caption{Several illustrations of the true solutions and their reconstructions in Example \ref{ex-d1e2}. (a) $f_\true^{\,\a}(x)=x$. (d) $f_\true^{\,\rd}(x)=1-|2x-1|$.}\label{fig-d1e2}
\end{figure}
\end{ex}

In the above examples, it is readily seen that even with quite coarse initial guesses $f_0$, the numerical reconstructions appear to be satisfactory in view of the ill-posedness of the inverse source problem. We evaluate the performance of our algorithm by analyzing the numerical results from the following aspects.

First, it is readily seen from Figures \ref{fig-d1e1} and \ref{fig-d1e2} that both of the above examples yield quite smooth reconstructions. In fact, according to the regularity result in Lemma \ref{lem-Lions}, one can expect an $H^1(\Om)$-regularity throughout the iteration \eqref{eq-itr-f}, as long as the initial guess $f_0$ and the known component $R(x,t)$ of the source term are sufficiently smooth. Such a smoothness, however, prevents us from proper identifications of non-smooth true solutions (see case (d) of Example \ref{ex-d1e2}).

Second, the reconstructed solutions appear more sensitive to the size of the observable subdomain $\om$ than to the data noise, but the non-monotonicity outside $\om$ is difficult to reconstruct. The influence of the smallness of $\om$ is witnessed from the second part of Table \ref{tab-d1e1}, which obviously comes from the limited information captured in $\om$. On the other hand, cases (a)--(c) in Example \ref{ex-d1e2} imply a tendency that the better the monotonicity of $f_\true$ is, the more accurate the identification will be, and the convergence will also be faster. In conclusion, although the conditional stability of the reconstruction is guaranteed by Lemma \ref{lem-inter}, in practice the signal strength from $\Om\setminus\ov\om$ is overwhelmed by that inside $\om$, so that the behavior of $f$ outside $\om$ cannot remarkably influence the observation data in $\om\times(0,T)$.

Third, Example \ref{ex-d1e1} suggests a considerably strong robustness of our algorithm against the measurement error. Actually, one can see from the first part of Table \ref{tab-d1e1} that the relative errors of the reconstructions only increase temperately as the observation noises are doubled. This phenomenon can be explained as follows. Suppose that the $m$-th iteration $f_m$ is of certain regularity, say $f_m\in L^2(\Om)$. According to Lemma \ref{lem-Lions}, the solution $u(f_m)$ to \eqref{eq-gov-w} should be sufficiently smooth, namely $u(f_m)\in\cD$ (see \eqref{eq-def-D} for the definition of $\cD$). Since the iteration \eqref{eq-itr-f} in principle aims at minimizing the surrogate functional $J^s(\,\cdot\,,f_m)$, it turns out that $u(f_m)|_{\om\times(0,T)}$ tends to take an averaged state of $u^\de$ in a sense that the error can be minimized. Therefore, provided that the observation data keep oscillating around the accurate ones, the reconstruction performs stably and insensitively in spite of the noise amplitude to a certain extent.

\Subsection{Two-dimensional examples}

Now we turn to the case of $n=2$. Without lose of generality, we always generate the subdomain $\om$ by removing a closed rectangle in $\Om=(0,1)^2$ whose edges are parallel to the coordinate axes. Due to the geometry condition for the reconstruction, $\ov\om$ should include at least two adjacent edges of $\ov\Om\,$. Simultaneously, the condition $T>\diam(\Om\setminus\ov\om)$ implies that the time duration $T$ should be longer than the diagonal of the removed rectangle. In the sequel, the largest size of such rectangles will be taken as $0.9^2$, and hence we will set $T=1.3>0.9\times\sqrt2$ in all tests for consistency. As before, we set the step size as $0.01$ and divide the space-time region as a $101^2\times131$ mesh in computation. Since the one-dimensional examples suggest that the reconstruction is insensitive to the noise, the noise level is always set as $5\%$ of the amplitude of $u(f_\true)$.

\begin{ex}\label{ex-d2e1}
In the first two-dimensional example, we fix $R(x,t)=5+\pi^2t^2$ and
\[f_\true(x)=f_\true(x_1,x_2)=\f12\cos(\pi x_1)\cos(\pi x_2)+1.\]
We test the algorithm by changing the subdomain $\om$ as follows. First, we keep the coverage $\pa\Om\subset\ov\om$ and reduce its thickness from $0.2$, $0.1$ to $0.05$, that is, we take $\om=\Om\setminus[0.2,0.8]^2$, $\om=\Om\setminus[0.1,0.9]^2$ and $\om=\Om\setminus[0.05,0.95]^2$ subsequently. Next, we fix the thickness as $0.1$ and reduce the coverage of $\pa\Om$ from $3$ edges to $2$, for instance, we choose $\om=\Om\setminus[0.1,1]\times[0.1,0.9]$ and $\om=\Om\setminus[0.1,1]^2$. The choices of $\om$ and various parameters as well as the corresponding numerical performances are listed in Table \ref{tab-d2e1}. The surface plots of several representative reconstructions $f_M$ are illustrated in Figure \ref{fig-d2e1}.
\begin{table}[htbp]\centering
\caption{Parameters and corresponding numerical performances in Example \ref{ex-d2e1} under various choices of observable subdomains.}\label{tab-d2e1}
\begin{tabular}{ccc|ccc|c}
\hline\hline
$\om$ & $K$ & $f_0$ & $M$ & $\mathrm{err}$ & elapsed time (s) & illustration\\
\hline
$\Om\setminus[0.2,0.8]^2$ & $3$ & $1$ & $31$ & $0.98\%$ & $30.28$ &\\
$\Om\setminus[0.1,0.9]^2$ & $1.7$ & $1$ & $28$ & $2.29\%$ & $28.05$ & Figure \ref{fig-d2e1}(a)\\
$\Om\setminus[0.05,0.95]^2$ & $1$ & $1$ & $27$ & $2.96\%$ & $26.35$\\
$\Om\setminus[0.1,1]\times[0.1,0.9]$ & $1.3$ & $1$ & $27$ & $3.46\%$ & $26.40$\\
$\Om\setminus[0.1,1]^2$ & $1$ & $1.5$ & $74$ & $7.53\%$ & $71.58$ & Figure \ref{fig-d2e1}(b)\\
\hline\hline
\end{tabular}
\end{table}

\begin{figure}[htbp]\centering
\includegraphics[trim=5mm 5mm 10mm 0mm,clip=true,width=.45\textwidth]{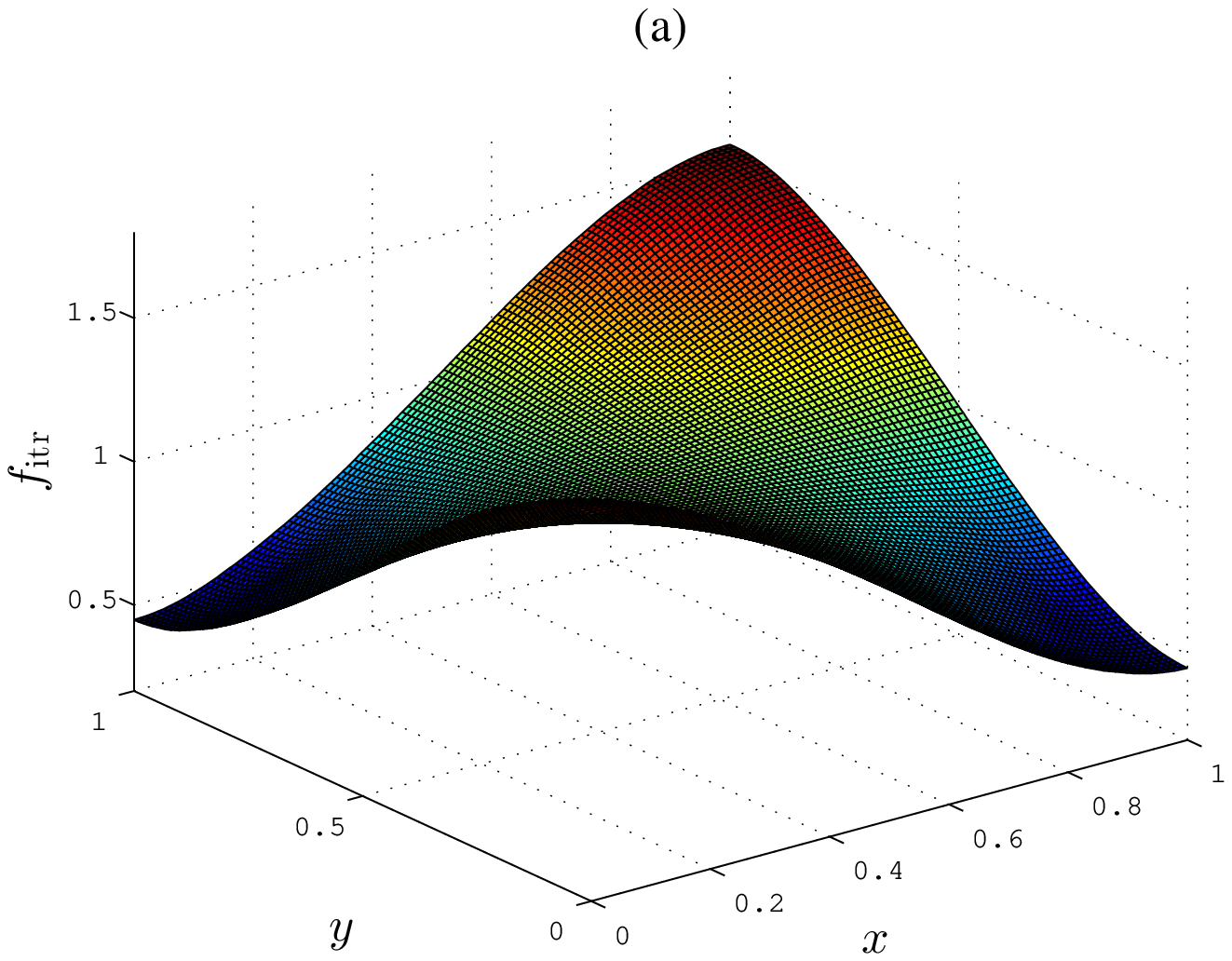}\qquad\includegraphics[trim=5mm 5mm 10mm 0mm,clip=true,width=.45\textwidth]{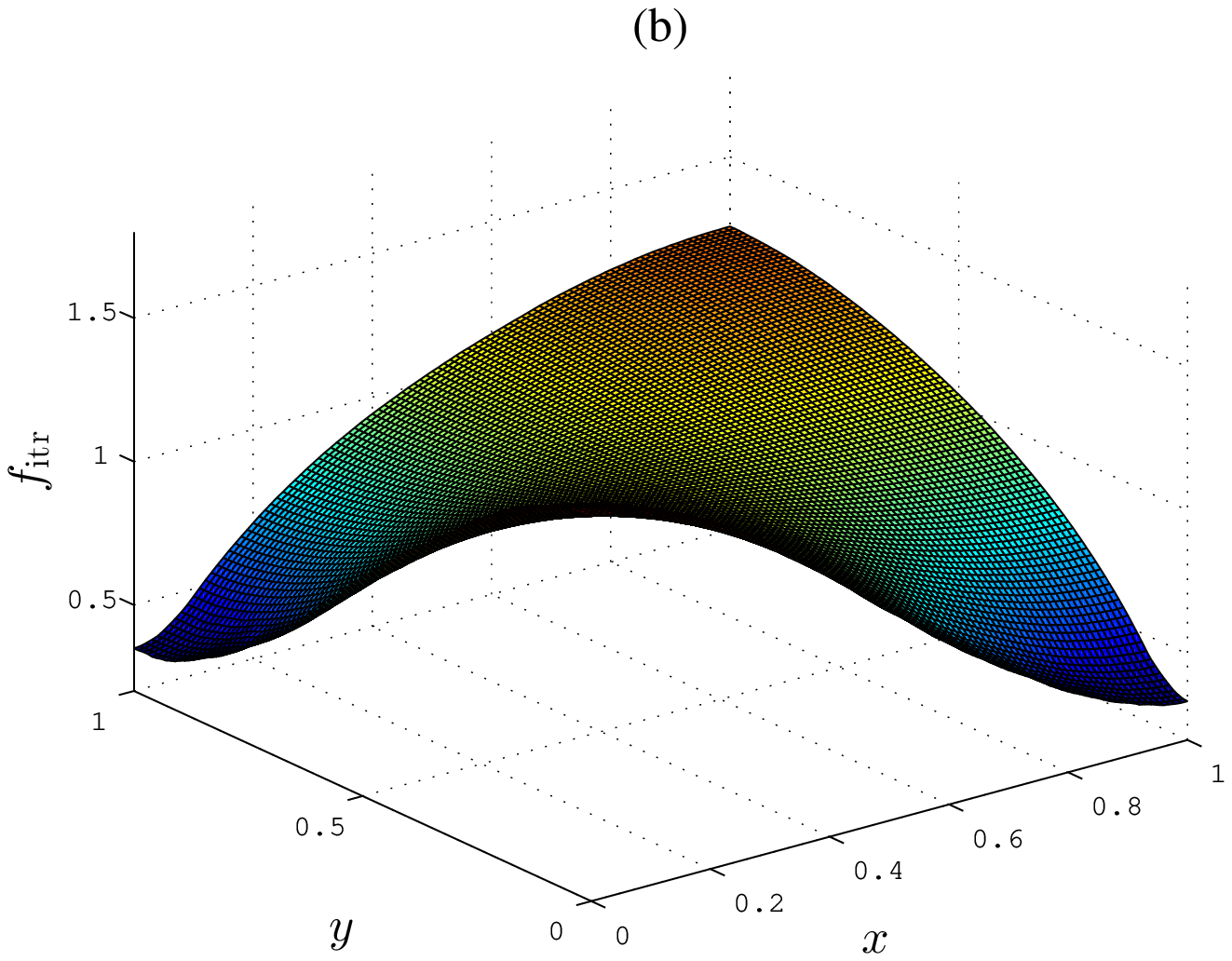}\\
\caption{Surface plots of several representative reconstructed solutions in Example \ref{ex-d2e1} with different choices of the observable subdomain $\om$. (a) $\om=\Om\setminus[0.1,0.9]^2$. (b) $\om=\Om\setminus[0.1,1]^2$.}\label{fig-d2e1}
\end{figure}
\end{ex}

\begin{ex}\label{ex-d2e2}
Parallelly to Example \ref{ex-d1e2} for the one-dimensional case, we investigate the influence of the monotonicity of $f_\true$ upon the numerical performance. To this end, we fix $R(x,t)=R(x_1,x_2,t)=x_1-x_2+3t+2$ and select three true solutions
\begin{align}
f_\true^{\,\a}(x) & =f_\true^{\,\a}(x_1,x_2)=\f12\cos(\pi x_1)+1,\label{eq-ftrue-d2a}\\
f_\true^{\,\rb}(x) & =f_\true^{\,\rb}(x_1,x_2)=3-\exp\left(1-\f{x_1+x_2}2\right),\label{eq-ftrue-d2b}\\
f_\true^{\,\rc}(x) & =f_\true^{\,\rc}(x_1,x_2)=\f12\cos(\pi x_1)\cos(2\pi x_2)+1.\label{eq-ftrue-d2c}
\end{align}
Here we take $\om=\Om\setminus[0.1,0.9]\times[0,0.9]$ as an intermediate choice, and set $K=0.27$, $f_0\equiv1$. The numbers $M$ of iterations and relative errors are shown in Table \ref{tab-d2e2}. As a typical example, we show the surface plots of $f_\true^{\,\rc}$ and its reconstruction in Figure \ref{fig-d2e2}.
\begin{table}[htbp]\centering
\caption{Numerical performances of the reconstructions in Example \ref{ex-d2e2} for various choices of true solutions.}\label{tab-d2e2}
\begin{tabular}{c|cc|c}
\hline\hline
$f_\true$ & $M$ & $\mathrm{err}$ & illustration\\
\hline
$f_\true^{\,\a}$ (see \eqref{eq-ftrue-d2a}) & $33$ & $2.70\%$\\
$f_\true^{\,\rb}$ (see \eqref{eq-ftrue-d2b}) & $41$ & $2.97\%$\\
$f_\true^{\,\rc}$ (see \eqref{eq-ftrue-d2c}) & $119$ & $7.22\%$ & Figure \ref{fig-d2e2}\\
\hline\hline
\end{tabular}
\end{table}

\begin{figure}[htbp]\centering
\includegraphics[trim=5mm 5mm 10mm 8mm,clip=true,width=.45\textwidth]{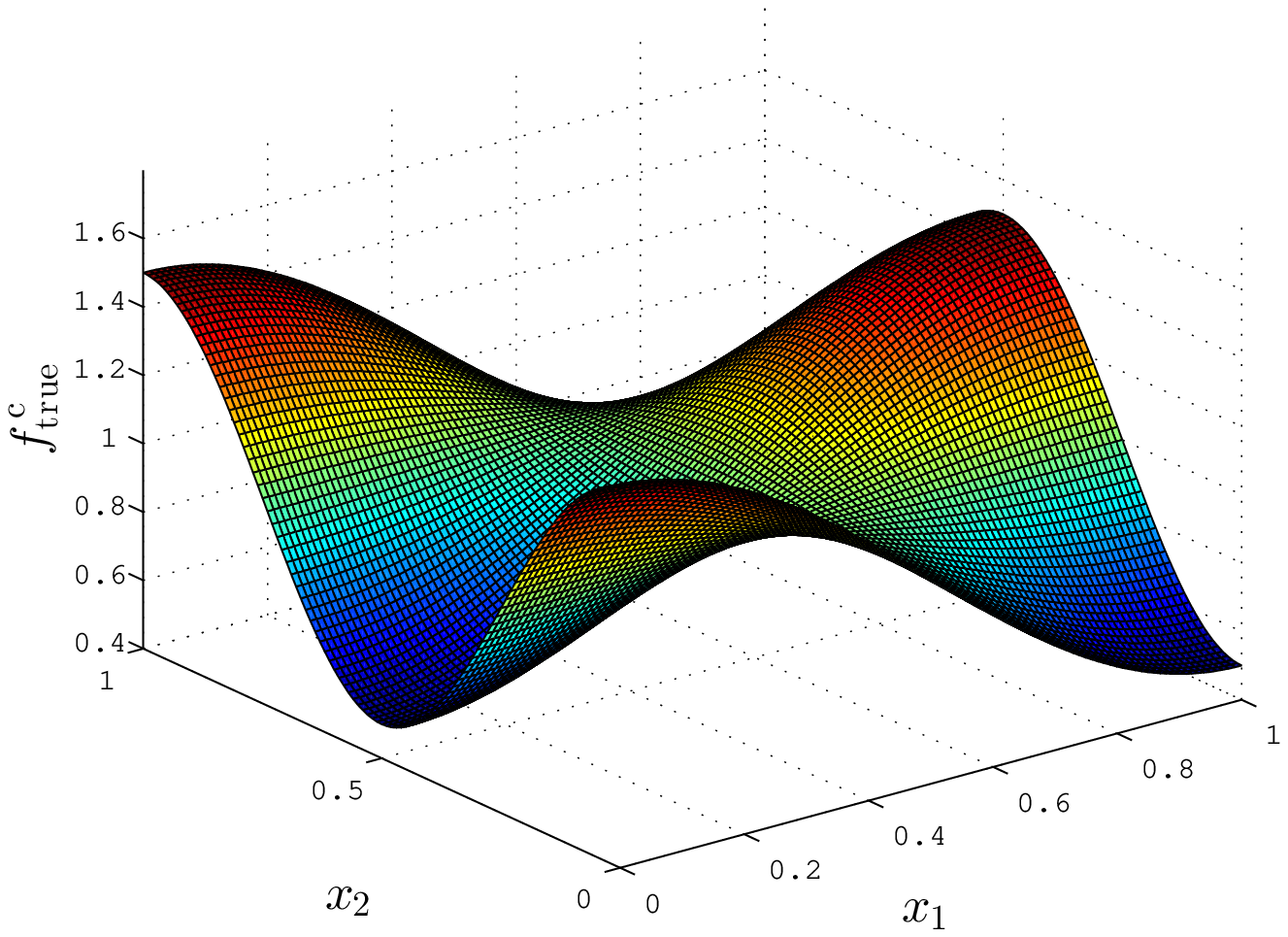}\qquad\includegraphics[trim=5mm 5mm 10mm 8mm,clip=true,width=.45\textwidth]{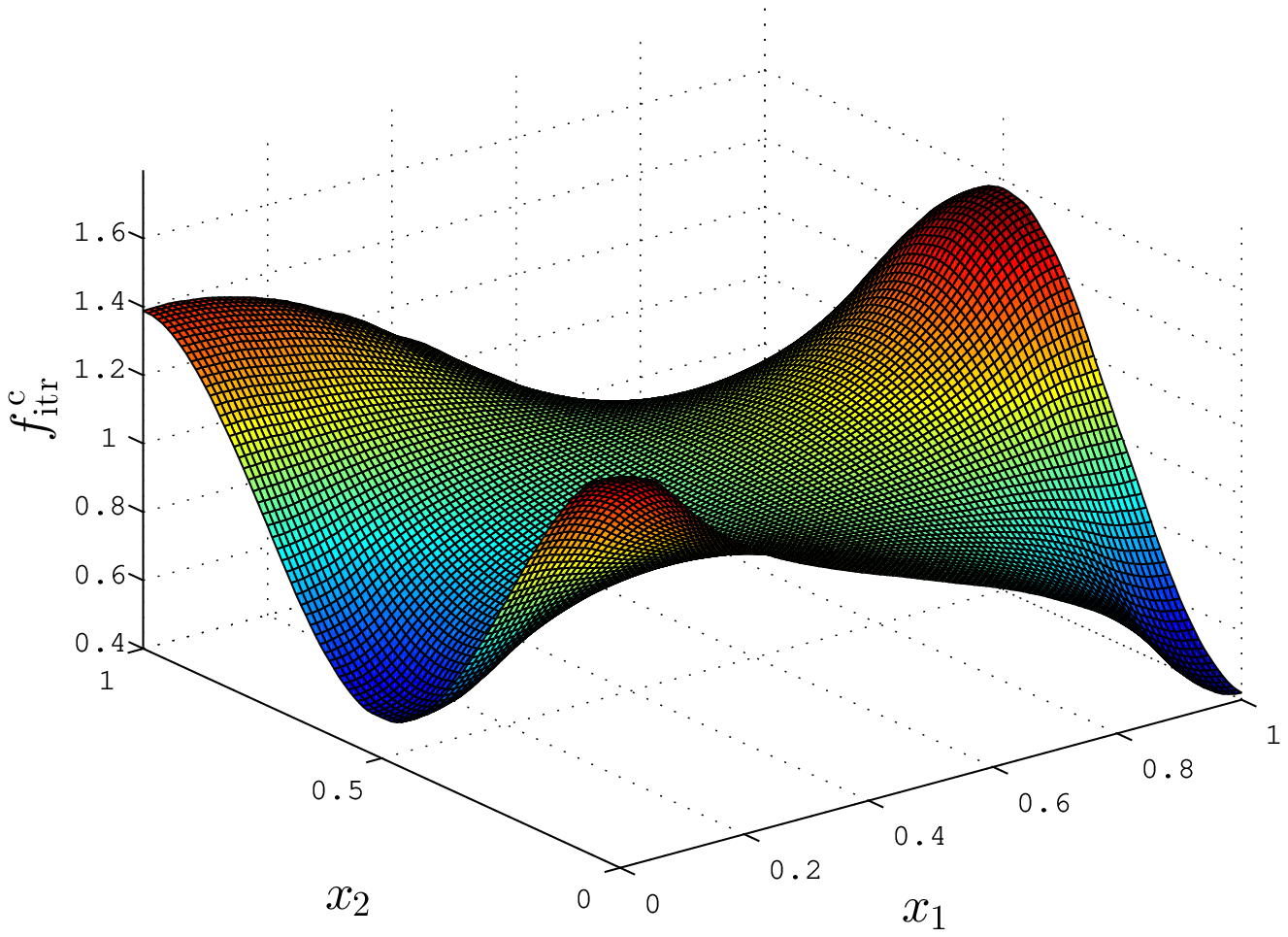}\\
\caption{Surface plots of $f_\true^{\,\rc}$ (left) and its reconstruction (right) in Example \ref{ex-d2e2}.}\label{fig-d2e2}
\end{figure}
\end{ex}

As expected, the above two-dimensional examples inherit mostly those phenomena observed in their one-dimensional counterparts and here we will not repeat the discussion again. Nevertheless we shall mention that, other than the thickness of $\om$, the iteration steps and relative errors of the numerical reconstructions also depend heavily on how much $\ov\om$ can cover $\pa\Om$. Although in all tests the algorithm performs quite well when $\pa\Om$ is included in $\ov\om$, it performs worse when $\ov\om$ only covers $3$ or $2$ edges of $\pa\Om$ in the sense of larger relative errors and more iteration steps. We can see clearly from Figure \ref{fig-d2e1}(b) and Figure \ref{fig-d2e2} that the numerical solutions $f_M$ fail to match with $f_\true$ especially on the uncovered edge, though $f_\true$ are well-reconstructed inside $\om$. In particular, Figure \ref{fig-d2e1}(b) indicates a dramatic difference at the uncovered corner $x=(1,1)$ when $\om=\Om\setminus[0.1,1]^2$. These demonstrate again the numerical ill-posedness regardless of the fact that the theoretical stability (Lemma \ref{lem-inter}) is valid under condition \eqref{eq-cond}.

\Subsection{Three-dimensional examples}

Finally, we proceed to the three-dimensional reconstruction of the source term $f_\true$. Similarly to the previous subsection, the subdomain $\om$ is generated by removing a closed cube in $\Om=(0,1)^3$ whose edges are parallel to the coordinate axes. According to the geometry condition for the reconstruction, $\ov\om$ should include at least three mutually adjacent faces of $\ov\Om\,$. Since the largest size of the removed cubes will be $0.96^3$, we set $T=1.7>0.96\times\sqrt3$ in all numerical tests in order to guarantee the condition $T>\diam(\Om\setminus\ov\om)$. Considering the computational complexity for $n=3$, we enlarge the mesh size in space and time as $0.02$ to produce a $51^3\times86$ mesh for $\ov\Om\times[0,T]$. As before, we still set the noise level as $5\%$ of the amplitude of $u(f_\true)$.

\begin{ex}\label{ex-d3e1}
Fix $R(x,t)=2+3\pi^2t^2$ and
\[f_\true(x)=f_\true(x_1,x_2,x_3)=\f12\cos(\pi x_1)\cos(\pi x_2)\cos(\pi x_3)+1.\]
Similarly to the previous subsections, we study the influence of the choice of $\om$ upon the numerical performance. First we keep the coverage $\pa\Om\subset\ov\om$ and use the thicknesses $0.08$ and $0.04$ to generate $\om=\Om\setminus[0.08,0.92]^3$ and $\om=\Om\setminus[0.04,0.96]^3$. Next we fix a thickness of $0.04$ and reduce the faces of $\pa\Om$ that $\ov\om$ covers from $5$ till $3$ which is the minimum possible, that is, e.g.\! $\om=\Om\setminus[0.04,0.96]^2\times[0.04,1]$, $\om=\Om\setminus[0.04,0.96]\times[0.04,1]^2$ and $\om=\Om\setminus[0.04,1]^3$. The choices of $\om$ and various parameters as well as the corresponding numerical performances are listed in Table \ref{tab-d3e1}.
\begin{table}[htbp]\centering
\caption{Parameters and corresponding numerical performances in Example \ref{ex-d3e1} under various choices of observable subdomains.}\label{tab-d3e1}
\begin{tabular}{ccc|ccc}
\hline\hline
$\om$ & $K$ & $f_0$ & $M$ & $\mathrm{err}$ & elapsed time (s)\\
\hline
$\Om\setminus[0.08,0.92]^3$ & $22$ & $1$ & $40$ & $1.83\%$ & $608.68$\\
$\Om\setminus[0.04,0.96]^3$ & $12$ & $1$ & $38$ & $2.49\%$ & $496.60$\\
$\Om\setminus[0.04,0.96]^2\times[0.04,1]$ & $10$ & $1$ & $38$ & $3.07\%$ & $546.37$\\
$\Om\setminus[0.04,0.96]\times[0.04,1]^2$ & $7.5$ & $1$ & $40$ & $3.87\%$ & $520.20$\\
$\Om\setminus[0.04,1]^3$ & $6$ & $0.5$ & $39$ & $8.84\%$ & $515.18$\\
\hline\hline
\end{tabular}
\end{table}
\end{ex}

\begin{ex}\label{ex-d3e2}
Finally, as before we test our algorithm by selecting true solutions with different degrees of monotonicity. We fix $R(x,t)=5+\pi^2t^2$, and choose
\begin{align}
f_\true^{\,\a}(x) & =f_\true^{\,\a}(x_1,x_2,x_3)=\left(x_1-\f15\right)\left(x_2-\f12\right)^2-\f12\cos(\pi x_1)\,x_3+\f12x_2\,\e^{-x_3},\label{eq-ftrue-d3a}\\
f_\true^{\,\rb}(x) & =f_\true^{\,\rb}(x_1,x_2,x_3)=\f12\cos(\pi x_1)\cos(2\pi x_2)\cos(\pi x_3)+1,\label{eq-ftrue-d3b}\\
f_\true^{\,\rc}(x) & =f_\true^{\,\rc}(x_1,x_2,x_3)=\f12\cos(\pi x_1)\cos(2\pi x_2)\cos(2\pi x_3)+1.\label{eq-ftrue-d3c}
\end{align}
Here we set $\om=\Om\setminus[0.04,0.96]\times[0.04,1]^2$, $K=3.5$ and $f_0\equiv1$. The numerical performances are listed Table \ref{tab-d3e2}.
\begin{table}[htbp]\centering
\caption{Numerical performances of the reconstructions in Example \ref{ex-d3e2} for various choices of true solutions.}\label{tab-d3e2}
\begin{tabular}{c|cc}
\hline\hline
$f_\true$ & $M$ & $\mathrm{err}$\\
\hline
$f_\true^{\,\a}$ (see \eqref{eq-ftrue-d3a}) & $64$ & $3.45\%$\\
$f_\true^{\,\rb}$ (see \eqref{eq-ftrue-d3b}) & $120$ & $7.88\%$\\
$f_\true^{\,\rc}$ (see \eqref{eq-ftrue-d3c}) & $101$ & $11.98\%$\\
\hline\hline
\end{tabular}
\end{table}
\end{ex}

Again, the three-dimensional examples show almost identical behaviors to that in lower dimensional cases. In summary, thinner observable subdomains $\om$ result in worse reconstructions, and their coverage of $\pa\Om$ also dominates the numerical performance to a great extent. On the other hand, the oscillation of $f_\true$ in $\Om\setminus\ov\om$ is extremely difficult to recover.

\Section{Concluding remarks}\label{sec-con}

Motivated by the global Lipschitz stability obtained in \cite{IY01a,IY01b} for the inverse source problem for wave equations, in this paper we investigate the same problem for the general hyperbolic equation with time-dependent principal part. Establishing a new Carleman estimate which dominates the solution up to second order derivatives, we prove the local H\"older stability for both cases of partial boundary observation and partial interior observation under suitable assumptions on the observation time and observable subboundary or subdomain.

Numerically, we consider the inverse source problem for a wave equation with partial interior measurements. Reformulating the problem as a minimization problem, we find the variational equation for the minimizer by using the corresponding adjoint system, which leads us to the iterative thresholding algorithm. Substantial numerical tests illustrate that the proposed algorithm is very accurate and efficient. In particular, our algorithm is considerably robust against the measurement error, but is sensitive to the size of the observable domain.

As a related work, in \cite{LJY15} we considered the same type of inverse source problem for the double wave equation
\[(\pa_t^2-\tri)^2u(x,t)=f(x)R(x,t)\]
which describes the three-dimensional time cone model. Similarly to the present paper, we also established a Lipschitz stability result and develop a corresponding iteration thresholding algorithm in that case. Our future work includes the application of the proposed iterative thresholding algorithm to inverse source problems for various kinds of partial difference equations with different types of observation data. For example, we will study similar problems for wave equations with final observation and fractional diffusion equations with partial interior observation. At the same time, we will also discuss the theoretical uniqueness and stability.

\appendix
\Section{Proof of Lemma \ref{lem-Carleman}}\label{sec-app}

Now we establish the key Carleman estimates \eqref{eq-CE-1}--\eqref{eq-CE-2} for the hyperbolic operator $\pa_t^2-\cA(t)$ on the basis of the existing Carleman estimate \eqref{eq-CE-0}.

Thanks to the large parameter $s>0$, it is sufficient to show Lemma \ref{lem-Carleman} in the case of $b=c=0$, that is,
\[F=(\pa_t^2-\cA(t))u=\pa_t^2u-\rdiv(a\nb u)\quad\mbox{in }Q,\]
Carleman estimate \eqref{eq-CE-0} yields
\begin{equation}\label{eq-CE-0'}
\int_{Q_\de}s\left(|\pa_tu|^2+|\nb u|^2+s^2|u|^2\right)\e^{2s\vp}\,\rd x\rd t\le C\int_{Q_\de}|F|^2\,\e^{2s\vp}\,\rd x\rd t.
\end{equation}
Setting $u_1:=\pa_tu$ and recalling the notation $\cA'(t)$ introduced in \eqref{eq-def-A'}, we immediately see
\[(\pa_t^2-\cA(t))u_1=\pa_tF+\cA'(t)u,\]
Applying \eqref{eq-CE-0} to $u_1$ and noting $|\cA'(t)u|^2\le C\sum_{i,j=1}^n|\pa_i\pa_ju|^2$, we have
\begin{align}
\int_{Q_\de}|\pa_t^2u|^2\,\e^{2s\vp}\,\rd x\rd t & =\int_{Q_\de}|\pa_tu_1|^2\,\e^{2s\vp}\,\rd x\rd t\le\f Cs\int_{Q_\de}\left|\pa_tF+\cA'(t)u\right|^2\e^{2s\vp}\,\rd x\rd t\nonumber\\
& \le C\int_{Q_\de}\f1s|\pa_tF|^2\,\e^{2s\vp}\,\rd x\rd t+\f Cs\int_{Q_\de}\sum_{i,j=1}^n|\pa_i\pa_ju|^2\,\e^{2s\vp}\,\rd x\rd t.\label{eq-est-utt}
\end{align}
Since $\cA(t)u=\pa_t^2u-F$ in $Q$,
we have
\begin{align}
& \quad\,\,\int_{Q_\de}|\cA(t)u|^2\,\e^{2s\vp}\,\rd x\rd t\le2\int_{Q_\de}\left(|\pa_t^2u|^2+|F|^2\right)\e^{2s\vp}\,\rd x\rd t\nonumber\\
& \le C\int_{Q_\de}\left(\f1s|\pa_tF|^2+|F|^2\right)\e^{2s\vp}\,\rd x\rd t+\f Cs\int_{Q_\de}\sum_{i,j=1}^n|\pa_i\pa_ju|^2\,\e^{2s\vp}\,\rd x\rd t.\label{eq-est-Au}
\end{align}
On the other hand, we calculate
\begin{align}
\pa_j(u\,\e^{s\vp}) & =\{\pa_ju+s(\pa_j\vp)u\}\,\e^{s\vp},\nonumber\\
\pa_i\pa_j(u\,\e^{s\vp}) & =\{\pa_i\pa_ju+s((\pa_i\vp)\pa_ju+(\pa_j\vp)\pa_iu+(\pa_i\pa_j\vp)u)+s^2(\pa_i\vp)(\pa_j\vp)u\}\,\e^{s\vp},\label{eq-uij}
\end{align}
indicating
\[\cA(t)(u\,e^{s\vp})=\left\{\cA(t)u+2s(a\nb\vp\cdot\nb u)+s(\cA(t)\vp)u+s^2(a\nb\vp\cdot\nb\vp)u\right\}\e^{s\vp}\quad\mbox{in }Q.\]
Therefore, there holds
\[|\cA(t)(u\,\e^{s\vp})|^2\le C\left(|\cA(t)u|^2+s^2|\nb u|^2+s^4|u|^2\right)\e^{2s\vp}\quad\mbox{in }Q\]
and thus
\[\int_\Om|\cA(t)(u\,\e^{s\vp})|^2\,\rd x\le C\int_\Om\left(|\cA(t)u|^2+s^2|\nb u|^2+s^4|u|^2\right)\e^{2s\vp}\,\rd x\quad(-T\le t\le T).\]
Applying the a priori estimate for the Dirichlet problem for the elliptic operator $\cA(t)$ (see, e.g., Gilbarg and Trudinger \cite{GT01}), we obtain for each fixed $t\in[-T,T]$ that
\begin{equation}\label{eq-est-uij}
\int_\Om\sum_{i,j=1}^n|\pa_i\pa_j(u\,\e^{s\vp})|^2\,\rd x\le C\int_\Om\left(|\cA(t)u|^2+s^2|\nb u|^2+s^4|u|^2\right)\e^{2s\vp}\,\rd x.
\end{equation}
Here, by $a\in W^{1,\infty}(-T,T;W^{1,\infty}(\Om))\subset C([-T,T];W^{1,\infty}(\Om))$, we note that we can choose $C>0$ uniformly in $t\in[-T,T]$. Then it follows from \eqref{eq-uij} and \eqref{eq-est-uij} that
\begin{align*}
\int_\Om\sum_{i,j=1}^n|\pa_i\pa_ju|^2\,\e^{2s\vp}\,\rd x & \le C\int_\Om\sum_{i,j=1}^n|\pa_i\pa_j(u\,\e^{s\vp})|^2\,\rd x+C\int_\Om\left(s^2|\nb u|^2+s^4|u|^2\right)\e^{2s\vp}\,\rd x\\
& \le C\int_\Om\left(|\cA(t)u|^2+s^2|\nb u|^2+s^4|u|^2\right)\e^{2s\vp}\,\rd x\quad(-T\le t\le T).
\end{align*}
Integrating over $(-T,T)$ and noting $\supp\,u\subset Q_\de$, we have
\begin{align*}
& \quad\,\,\int_{Q_\de}\sum_{i,j=1}^n|\pa_i\pa_ju|^2\,\e^{2s\vp}\,\rd x\rd t\le C\int_{Q_\de}\left(|\cA(t)u|^2+s^2|\nb u|^2+s^4|u|^2\right)\e^{2s\vp}\,\rd x\rd t\\
& \le C\int_{Q_\de}\left(\f1s|\pa_tF|^2+|F|^2\right)\e^{2s\vp}\,\rd x\rd t+\f Cs\int_{Q_\de}\sum_{i,j=1}^n|\pa_i\pa_ju|^2\,\e^{2s\vp}\,\rd x\rd t+C\,s\int_{Q_\de}|F|^2\,\e^{2s\vp}\,\rd x\rd t\\
& \le C\int_{Q_\de}\left(\f1s|\pa_tF|^2+s\,|F|^2\right)\e^{2s\vp}\,\rd x\rd t+\f Cs\int_{Q_\de}\sum_{i,j=1}^n|\pa_i\pa_ju|^2\,\e^{2s\vp}\,\rd x\rd t,
\end{align*}
where the second inequality follows from \eqref{eq-est-Au} and \eqref{eq-CE-0'}. Choosing $s>0$ sufficiently large, we can absorb the last term on the right-hand side into the left-hand side, which implies \eqref{eq-CE-1}.

Finally, it is straightforward to obtain \eqref{eq-CE-2} by substituting \eqref{eq-CE-1} into \eqref{eq-est-utt}.

\end{document}